\documentclass{amsart}
\usepackage{epsfig}
\theoremstyle{plain}
\newtheorem{theorem}{Theorem}[section]
\newtheorem{lemma}[theorem]{Lemma}
\newtheorem{proposition}[theorem]{Proposition}
\newtheorem{corollary}[theorem]{Corollary}

\theoremstyle{definition}
\newtheorem{definition}[theorem]{Definition}
\newtheorem{remark}[theorem]{Remark}

\theoremstyle{remark}
\newtheorem*{remark*}{Remark}

 \def\bd{\partial}
 \def\C{{\mathbb C}}
 \def\R{{\mathbb R}}
 
 \def\H{{\mathbb H}}
 \def\Z{{\mathbb Z}}
 \def\g{{\mathfrak g}}

 \def\t{\theta}

\def\arcsinh{\mathop{\rm arcsinh}}

\def\arctanh{{\mathop{\rm arctanh}}}
\def\area{\mathop{\rm area}}

\def\vol{\mathop{\rm vol}}

\def\M{{X}} %name for cusped manifold 
\def\N{{M}} %name for cone-manifold   
\def\m{{m}} % meridian length on maximal tube boundary 
\def\h{{h}} % height of maximal tube boundary      
\def\b{{b}} %boundary term 
\def\l{{\ell}} % core geodesic length
\def\L{L} %name for geodesic length on cusp torus
\def\z{\zeta} 
\def\T{\tau}

\def\cal{\mathcal}
\def\Re{\mathop{\rm Re}}
\def\S{{S}}
\def\r{\rho}
\def\tw{{\rm tw}}

\title{Universal bounds for hyperbolic Dehn surgery}
\author{Craig D. Hodgson}
\address{Department of Mathematics and Statistics\\
University of Melbourne\\
Victoria 3010, Australia}
\email{cdh@ms.unimelb.edu.au}
\thanks{The research of the first author was partially
supported by grants from the ARC}
\author{Steven P. Kerckhoff}
\address{Department of Mathematics\\
Stanford University\\
Stanford, CA 94305, U.S.A.}
\email{spk@math.stanford.edu}
\thanks{The research of the second author was partially
supported by grants from the NSF}

\begin{document}

\begin{abstract}
This paper gives a quantitative version of Thurston's hyperbolic Dehn
surgery theorem. Applications include the first universal bounds on the 
number of non-hyperbolic Dehn fillings on a cusped hyperbolic 3-manifold,
and estimates on the changes in volume and core geodesic length during 
hyperbolic Dehn filling. The proofs involve the construction of a family 
of hyperbolic cone-manifold structures, using infinitesimal harmonic 
deformations and analysis of geometric limits.
\end{abstract}

\maketitle

\section{Introduction}\label{intro}

If $\M$ is a non-compact, finite volume orientable 
hyperbolic $3$-manifold, it
is the interior of a compact $3$-manifold with a finite number of
torus boundary components.  For each torus, there are an
infinite number of topologically distinct ways to attach a solid torus. Such
``Dehn fillings" are parameterized by pairs of relatively prime
integers, once a basis for the fundamental group of the torus
is chosen.  If each torus is filled, the resulting manifold is
closed.  A fundamental theorem of Thurston (\cite{thnotes})
states that, for all but a finite number of Dehn surgeries on
each boundary component, the resulting closed $3$-manifold has a
hyperbolic structure.  However, it was unknown whether or not
the number of such non-hyperbolic surgeries was bounded independent of
the original non-compact hyperbolic manifold.

In this paper we obtain a
universal upper bound on the number of non-hyperbolic
Dehn surgeries per boundary torus, independent of the
manifold $\M$. There are at most $60$ non-hyperbolic 
Dehn surgeries if there is only one cusp; if there are
multiple cusps, at most $114$ surgery curves must be excluded
from each boundary torus.  

These results should be compared with the known bounds on
the number of Dehn surgeries which yield 
manifolds which fail to be either irreducible or atoroidal or fail to have 
infinite fundamental group.  These are all necessary conditions for 
a $3$-manifold to be hyperbolic.  The hyperbolic geometry part of 
Thurston's geometrization conjecture states that these conditions 
should also be sufficient; i.e., that the interior of a compact, 
orientable $3$-manifold has a complete hyperbolic structure if and 
only if it is irreducible, atoroidal, and has infinite fundamental group.

It follows from the work of Gromov-Thurston
(\cite{GT}, see also \cite{BH}) that all but a universal number of 
surgeries on each torus yield $3$-manifolds which admit 
negatively curved metrics.  
More recent work by Lackenby \cite{La} and, independently, by
Agol \cite{agol}, similarly shows that for all but a universally
bounded number of surgeries on each torus the resulting manifolds
are irreducible with infinite word hyperbolic fundamental group.
Similar types of bounds using techniques
less comparable to those in this paper have been obtained by
Gordon, Luecke, Wu, Culler, Shalen, Boyer, Zhang and many others.
(See, for example, \cite{CGLS}, \cite{BZ} and the survey articles
 \cite{Go1}, \cite{Go2}.)
Negatively curved $3$-manifolds are irreducible, atoroidal
and have infinite fundamental groups. If the geometrization
conjecture were known to be true, it would imply that these manifolds
actually have hyperbolic metrics.  The same is true for irreducible
$3$-manifolds with infinite word hyperbolic fundamental group.
Thus, the above results would provide a universal bound on the
number of non-hyperbolic Dehn fillings.  However, without first establishing
the geometrization conjecture, no such conclusion is possible and other methods
are required.

The bound on the number of Dehn surgeries that fail to be negatively
curved comes from what is usually referred to as the ``$2\pi$-theorem".
It can be stated as follows:  Given a cusp in a complete, orientable
hyperbolic $3$-manifold $\M$, remove a horoball neighborhood of the
cusp, leaving a manifold with a boundary torus which has a
flat metric.  Let $\gamma$ be an isotopy class of simple closed
curve on this torus and let $\M(\gamma)$ denote $\M$ filled in
so that $\gamma$ bounds a disk.  Then the $2\pi$-theorem states
that, if the flat geodesic length of $\gamma$ on the torus is
greater than $2\pi$, then $\M(\gamma)$ can be given a metric
of negative curvature which agrees with the hyperbolic metric in 
the region outside the horoball.  The bound then follows from the 
fact that it is always possible to find an embedded horoball 
neighborhood with boundary torus whose shortest geodesic has 
length at least 1.  On such a torus there are a bounded number of 
isotopy classes of geodesics with length less than or equal to $2\pi$.

Similarly, Lackenby and Agol show that, if the flat geodesic length
is greater than $6$, then the Dehn filled manifold is irreducible with
infinite word hyperbolic fundamental group.  
Agol then uses the recent work of Cao-Meyerhoff
(\cite {cao-meyerhoff}), which provides an improved lower bound on
the area of the maximal embedded horotorus, to conclude that, when
there is a single cusp, at most $12$ surgeries fail to be irreducible
or infinite word hyperbolic.  This is remarkably close to the the largest known 
number of non-hyperbolic Dehn surgeries which is 10, occurring for
the complement of the figure-8 knot.

Our criterion for those surgery curves whose corresponding filled 
manifold is guaranteed to be hyperbolic is similar.  We consider the
{\em normalized length} of curves on the torus, 
measured after rescaling the metric on the torus to have area 1, i.e. 
${\rm normalized~length} = ({\rm geodesic~length})/\sqrt{{\rm torus~area}}$.
Our main result shows that, if the {normalized length} 
of $\gamma$ on the torus is sufficiently long, then it is
possible to deform the complete hyperbolic structure through
cone-manifold structures on $\M(\gamma)$ with $\gamma$ bounding a 
singular meridian disk until the cone angle reaches $2 \pi$.  This
gives a smooth hyperbolic structure on $\M(\gamma)$.
The important point here is that ``sufficiently long" is
universal, independent of $\M$.  As before, it is straightforward
to show that all but a universal number of isotopy classes of simple 
closed curves satisfy this normalized length condition.  

The condition in this case that the normalized length, rather than 
just the flat geodesic length, be long is probably not necessary,
but is an artifact of the proof.  

We will now give a rough outline of the proof.

We begin with a non-compact, finite volume hyperbolic $3$-manifold
$\M$, which, for simplicity, we assume has a single cusp.  In the general
case the cusps are handled independently.  The manifold $\M$ is the 
interior of a compact manifold which has a single torus boundary.  
Choose a simple closed curve $\gamma$ on the torus. We wish to put 
a hyperbolic structure on the closed manifold $\M(\gamma)$ obtained 
by Dehn filling.  The metric on the open manifold $\M$ is deformed 
through incomplete metrics whose metric completion is a singular metric 
on $\M(\gamma)$, called a cone metric.  (See \cite{HK1} for a 
detailed description of such metrics.)  The singular set is a simple 
closed geodesic at the core of the added solid torus.  For any plane 
orthogonal to this geodesic the disks of small radius around the 
intersection with the geodesic have the metric of a 2-dimensional 
hyperbolic cone with angle $\alpha$.  The angle $\alpha$ is the same 
at every point along the singular geodesic $\Sigma$ and is called 
the cone angle at $\Sigma$. The complete structure can be considered 
as a cone-manifold with angle $0$.  The cone angle is increased 
monotonically, and, if the angle of $2 \pi$ is reached, this defines 
a smooth hyperbolic metric on $\M(\gamma)$.

The theory developed in \cite{HK1} shows that it is always possible to
change the cone angle a small amount, either increase it or decrease
it.  Furthermore, this can be done in a unique way, at least locally.
The cone angles locally parameterize the set of cone-manifold structures
on $\M(\gamma)$.  In particular, there are no variations of the hyperbolic
metric which leave the cone angle fixed.  This property is referred
to as {\em local rigidity rel cone angles}. Thus, to choose a 
1-parameter family of cone angles is to choose a well-defined 
family of singular hyperbolic metrics on $\M(\gamma)$ of this type. 

Although there are always local variations of the cone-manifold structure,
the structure may degenerate in various ways as a family of
angles reaches a limit.  In order to find a smooth hyperbolic metric
on $\M(\gamma)$ it is necessary to show that no degeneration occurs
before the angle $2 \pi$ is attained.

The proof has two main parts, involving rather different types
of arguments.  One part is fairly analytic, showing that
under the normalized length hypothesis on $\gamma$, there is a lower bound 
to the tube radius for any of the cone-manifold 
structures on $\M(\gamma)$ with angle at most $2 \pi$.  
The second part consists of showing that, under certain geometric conditions,
most importantly the lower bound on the tube radius, no degeneration of the 
hyperbolic structure is possible.  This involves studying possible 
geometric limits where the tube radius condition restricts such limits to fairly 
tractable and well-understood types.

The argument showing that there is a lower bound to the tube
radius is based on the local rigidity theory for cone-manifolds
developed in \cite{HK1}.  Indeed, the key estimates are best viewed
as effective versions of local rigidity of cone-manifolds.
We choose a smooth parametrization of the increasing family of cone angles,
which uniquely determines a family of cone-manifold structures.
We then need to control the global behavior of these metrics.
The idea is first to form a model for the deformation in a neighborhood
of the singular locus which changes the cone angle in the prescribed fashion 
and then find estimates which bound
the deviation of the actual deformation from the model.

The main goal is to estimate the actual behavior of the holonomy of 
the fundamental group elements corresponding to the boundary torus.  
The holonomy representation of the meridian is simply an elliptic 
element which rotates by the cone angle so it suffices to understand 
the longitudinal holonomy.  We derive some estimates on the
complex length of the longitude in terms of the cone angle
which depend on the original geometry of the horospherical
torus, including the length of the meridian on the torus.
These results may be of independent interest.

The estimates are derived by analyzing boundary terms in a
Weitzenb\"ock formula for the infinitesimal change of metric
which arises from differentiating our family of cone metrics.
This formula is the basis for local rigidity of hyperbolic metrics
in dimensions 3 and higher (\cite{cal}, \cite{weil}) and of hyperbolic
cone-manifolds in dimension 3 (\cite{HK1}). Our estimates ultimately provide
a bound on the derivative of the ratio of the cone angle to the
hyperbolic length of the singular core curve of the cone-manifold.
The bound depends on the tube radius.  On the other hand, a
geometric packing argument shows that the change in the tube
radius can be controlled when the product of the cone angle
and the core length is small.

Putting these results together, we arrive at differential 
inequalities which provide strong control on the change in the 
geometry of the maximal tube around the singular geodesic, including
the tube radius.  The value of the normalized flat length of
the surgery curve on the maximal cusp torus for the complete 
structure gives the initial condition for the ratio of the
cone angle to the core length.  (Note:  The {\it ratio}
of the cone angle to the core length approaches a finite,
non-zero value even though they individually approach zero 
at the complete structure.)

The conclusion is that, if the initial value of the ratio is
large, then it will remain large and the product of the cone angle
and the core length will remain small.  The packing argument then shows
that there will be a lower bound to the tube radius.

This gives a proof of the following Theorem:

\begin{theorem}\label{thm1}  Let $\M$ be a complete, finite volume 
orientable hyperbolic $3$-manifold with one cusp and 
let $T$ be a horospherical torus which is embedded as
a cross-section to the cusp. Let
$\gamma$ be a simple closed curve on $T$, $\M(\gamma)$ the Dehn
filling with $\gamma$ as meridian.  Let $\M_\alpha(\gamma)$ be 
a cone-manifold structure on $\M(\gamma)$ with cone angle $\alpha$
along the core, $\Sigma$, of the added solid torus, obtained by increasing
the angle from the complete structure.
If the normalized length of $\gamma$ on $T$ is at least $7.515$, 
then there is a positive lower bound to the tube radius around $\Sigma$
for all $2\pi \ge \alpha \ge 0$.
\end{theorem}

%\medskip

This theorem doesn't guarantee that cone angle $2 \pi$ can actually
be reached, just that there 
is a lower bound to the tube radius over all
angles less than or equal to $2 \pi$ that {\em are} attained.  
That $2 \pi$ can actually be attained follows from the next theorem.

\begin{theorem} \label{thm2}  Let $\N_t$, $t\in [0,t_\infty)$, be
a smooth path of closed hyperbolic cone-manifold structures
on $(\N, \Sigma)$ with cone angle $\alpha_t$ along the singular locus $\Sigma$.
Suppose $\alpha_t \to \alpha \geq 0$ as $t \to t_\infty$, that the volumes of
the $\N_t$ are bounded above by $V_0$, and that there is a positive constant
$R_0$ such that there is an embedded tube of radius at least $R_0$ around
$\Sigma$ for all $t$.  Then the path extends continuously to $t=t_\infty$
so that as $t \to t_\infty$, $\N_t$ converges in the bilipschitz topology
to a cone-manifold structure $\N_\infty$ on $\N$ with cone angles
$\alpha$ along $\Sigma$.
\end{theorem}

Given $\M$ and $T$ as in Theorem \ref{thm1}, choose {\em any}
non-trivial simple closed curve $\gamma$ on $T$.  There is a
maximal sub-interval $J \subset [0, 2\pi]$ containing $0$ such
that there is a smooth family $\N_\alpha$, where $\alpha \in J$,
of hyperbolic cone-manifold structures on $\M(\gamma)$ with
cone angle $\alpha$.  Thurston's Dehn surgery theorem (\cite{thnotes})
implies that $J$ is non-empty and \cite [Theorem 4.8] {HK1} implies
that it is open.  Theorem \ref{thm2} implies that, with a lower bound
on the tube radii and an upper bound on the volume, the path of 
$\N_\alpha$'s can be extended continuously to the endpoint of $J$.  
Again, \cite [Theorem 4.8] {HK1}
implies that this extension can be made to be smooth.  Hence,
under these conditions $J$ will be closed.  By Schl\"afli's formula 
(\ref{schlafli}) the volumes decrease as the cone angles increase,
so they will clearly be bounded above.  Theorem \ref{thm1} provides
initial conditions on $\gamma$ which guarantee that there will be
a lower bound on the tube radii for all $\alpha \in J$.  Thus, assuming Theorems
\ref{thm1} and \ref{thm2}, we have proved:

\begin{theorem}  \label{thm3}
Let $\M$ be a complete, finite volume orientable hyperbolic $3$-manifold with 
one cusp, and let $T$ be a horospherical torus which is embedded as
a cross-section to the cusp of $\M$. Let
$\gamma$ be a simple closed curve on $T$ whose Euclidean geodesic length on
$T$ is denoted by $\L$.  If the normalized length of $\gamma$,  
$\hat \L = \displaystyle{\L \over \sqrt{{\rm area}(T)}}$, 
is at least $7.515$,
then the closed manifold $\M(\gamma)$ obtained by Dehn filling along $\gamma$ is
hyperbolic.
\end{theorem}

%\medskip
This result also gives a universal bound on the
number of non-hyperbolic Dehn fillings on a cusped hyperbolic $3$-manifold
$\M$, independent of $\M$.

\begin{corollary} \label{thm4} 
Let $\M$ be a complete, orientable hyperbolic $3$-manifold with 
one cusp. Then at most $60$ Dehn fillings on $\M$ yield manifolds which admit no
complete hyperbolic metric.
\end{corollary}

When there are multiple cusps the results (Theorem \ref{multiplecusp}) 
are only slightly weaker. 
Theorem \ref{thm2} holds without change. If there are $k$ cusps,
the cone angles $\alpha_t$ and $\alpha$ are simply interpreted as
$k$-tuples of angles.  Having tube radius at least $R$ is interpreted
as meaning that there are disjoint, embedded tubes of radius $R$ 
around all of components of the singular locus.  The conclusion of 
Theorem \ref{thm1} and hence of Theorem \ref{thm3} holds when there 
are multiple cusps as long as the normalized lengths of all of the 
meridian curves are at least $\sqrt{2} ~ 7.515 \approx 10.6273$. 
At most $114$ curves from each cusp need to be excluded.  In fact,
this can be refined to say that at most $60$ curves need to be 
excluded from one cusp and at most $114$ excluded from the remaining 
cusps.  The rest of the Dehn filled manifolds are hyperbolic. 

In the final section of the paper (Section \ref{compare}), we
prove that every closed hyperbolic manifold with a sufficiently short 
(length less than $.111$) closed geodesic can be obtained by the 
process studied in this paper.  Specifically,
if one removes a simple closed geodesic from a closed hyperbolic manifold,
the resulting manifold can be seen to have a complete, finite volume
hyperbolic structure.  We prove that, if the removed geodesic had
length less than $.111$, then the hyperbolic structure on the closed
manifold and that of the complement of the geodesic can be connected
by a smooth family of hyperbolic cone-manifolds, with angles varying monotonically
from $2 \pi$ to $0$.

Also in that section (Theorem \ref{volchange}), we prove inequalities bounding the
difference between the volume of a complete hyperbolic manifold
and certain closed hyperbolic manifolds obtained from it by
Dehn filling.  We see (Corollary \ref{volL}) that, for the
manifolds constructed in Theorem \ref{thm3}, this difference
is at most $0.329$.  Similarly, using known bounds on the
volume of cusped hyperbolic $3$-manifolds, we prove (Corollary
\ref{volshort}) that every closed $3$-manifold with a closed
geodesic of length at most $0.162$ has volume at least
$1.701$.

\medskip

This paper is organized as follows:

\medskip

In Section \ref{harmdef} we recall basic definitions for deformations
of hyperbolic structures and some necessary results from a previous
paper (\cite{HK1}).  We use these to derive our fundamental inequality
(Theorem \ref {keyboundtheorem}) for the variation of the length of the 
singular locus as the cone angle is changed.  Section \ref{geomlim}
analyzes the limiting behavior of sequences of hyperbolic cone-manifolds
under the hypothesis of a lower bound to the tube radius around the
singular locus.  The proof of Theorem \ref{thm2} is given in that section.
It is, for the most part, independent of the rest of the paper.  
In Section \ref{pack} we use a packing argument to
relate the tube radius to the length of the singular locus.
In Section \ref{tuberad} we combine this relation with the inequality
from Section \ref{harmdef} to derive initial conditions that ensure
that there will be a lower bound to the tube radius for all cone angles
between $0$ and $2\pi$.  In particular,
the proof of Theorem \ref{thm1} is completed in that section.

\section{Deformation models and changes in holonomy}\label{harmdef}

In this section we recall the description of an infinitesimal
change of hyperbolic structure in terms of bundle valued 1-forms
and the Weitzenb\"ock formula satisfied by such a form when it is
harmonic in a suitable sense.  We compute the boundary term
for this formula in some specific cases which will allow us
to estimate the infinitesimal changes in the holonomy representations
of peripheral elements of the fundamental group.

In order to discuss the analytic and geometric objects associated
to an infinitesimal deformation of a hyperbolic structure, we need
first to describe what we mean by a $1$-parameter family of hyperbolic
structures.

A hyperbolic structure on an $n$-manifold $\M$
is determined by local charts modelled on $\H^n$ whose overlap
maps are restrictions of global isometries of $\H^n$.  These
determine, via analytic continuation,
a map $\Phi :\tilde \M \to \H^n$ from the universal cover
$\tilde \M$ of $\M$ to $\H^n$, called the {\em developing map},
which is determined uniquely up to post-multiplication by an
element of $G = {\rm isom} (\H^n)$. The developing map satisfies
the equivariance property $\Phi(\gamma m) = \rho(\gamma)\Phi(m)$,
for all $m\in \tilde \M$, $\gamma \in \pi_1(\M)$,
where $\pi_1(\M)$ acts on $\tilde \M$ by covering transformations,
and $\rho: \pi_1(\M) \to G$
is the {\em holonomy representation} of the structure.
The developing map also determines the hyperbolic metric on $\tilde \M$
by pulling back the hyperbolic metric on $\H^n$.
(See \cite{thbook} and \cite{Rat} for a complete discussion of these ideas.)

We say that two hyperbolic structures are {\em equivalent}
if there is a diffeomorphism $f$ from $\M$ to itself taking one
structure to the other.  We will use the term ``hyperbolic
structure" to mean such an equivalence class.
A {\em 1-parameter family}, $\M_t$, of hyperbolic structures defines
a 1-parameter family of developing maps
$\Phi_t:\tilde \M \to \H^n$,
where two families are equivalent under the relation
$\Phi_t \equiv k_t \Phi_t \tilde f_t$
where $k_t$ are isometries of $\H^n$ and $\tilde f_t$ are
lifts of diffeomorphisms $f_t$ from $\M$ to itself.
We assume that $k_0$ and $\tilde f_0$ are the identity,
and write $\Phi_0=\Phi$. 
All of the maps here are assumed to be smooth and to vary
smoothly with respect to $t$.

The tangent vector to a smooth family of hyperbolic
structures will be called an {\em infinitesimal deformation}.
The derivative at $t = 0$ of a a 1-parameter family of developing
maps $\Phi_t:\tilde \M \to \H^n$ defines a map
$\dot \Phi: ~ \tilde \M \to T\H^n$.  For any point
$m \in \tilde \M$, $\Phi_t(m)$ is a curve
in $\H^n$ describing how the image of $m$ is moving under the
developing maps; $\dot \Phi (m)$ is the initial tangent vector to the
curve.

We will identify $\tilde \M$ locally with $\H^n$ and $T \tilde \M$
locally with $T\H^n$ via the initial developing map $\Phi$.
Note that this identification is generally not a homeomorphism unless
the hyperbolic structure is complete.  However, it is a local
diffeomorphism, providing identification of small open sets
in $\tilde \M$ with ones in $\H^n$.

In particular, each point $m \in \tilde \M$ has a neighborhood $U$ where
$\Psi_t = \Phi^{-1} \circ \Phi_t : U \to \tilde \M$ is defined,
and the derivative at $t=0$ defines a vector field on $\tilde \M$,
$v = \dot \Psi : \tilde \M \to T \tilde \M$.
This vector field
determines the variation in developing maps
since $\dot \Phi = d\Phi \circ v$, and also determines
the variation in metric as follows. Let $g_t$ be the
hyperbolic metric on $\tilde \M$ obtained by pulling back the
hyperbolic metric on $\H^n$ via $\Phi_t$ and put $g_0=g$. Then
$g_t = \Psi_t^* g$ and the variation in metrics
$\dot g = {d  g_t \over dt} |_{t=0}$ is
the Lie derivative, ${\cal L}_v g$, of the initial metric $g$ along $v$.

Covariant differentiation of
the vector field $v$ gives a $T \tilde \M$ valued
1-form on $\tilde \M$, $\nabla v : T\tilde \M \to T \tilde \M$,
defined by $\nabla v(x) = \nabla_x v$ for $x \in T\tilde \M$.
We can decompose $\nabla v$ at each point into a symmetric part
and skew-symmetric part.
The {\em symmetric part},
$\tilde\eta = (\nabla v)_{sym}$,
represents the infinitesimal change in metric, since
$$\dot g (x,y) = {\cal L}_v g (x,y) =
g(\nabla_x v,y) + g(x,\nabla_y v) = 2 g(\tilde \eta(x),y)$$
for $x,y \in T \tilde \M$. In particular, $\tilde\eta$
descends to a well-defined $T\M$-valued 1-form $\eta$ on $\M$.
The {\em skew-symmetric} part $(\nabla v)_{skew}$ is the
{\em curl} of the vector field $v$, and its
value at $m \in \tilde \M$ represents
the effect of an infinitesimal rotation about $m$.

To connect this discussion of infinitesimal deformations
with cohomology theory, we consider the
Lie algebra  $\g$ of $G = {\rm isom} (\H^n)$
as vector fields on $\H^n$ representing infinitesimal
isometries of $\H^n$. Pulling back these vector fields
via the initial developing map $\Phi$
gives locally defined infinitesimal
isometries on $\tilde \M$ and on $\M$.

Let $\tilde E, E$ denote the vector bundles over
$\tilde \M, \M$ respectively of (germs of)
infinitesimal isometries. Then we can regard $\tilde E$ as the product
bundle with total space $\tilde \M \times {\g}$,
and $E$ is isomorphic to $(\tilde \M \times {\g})/\!\!\sim$
where $(m,v) \sim (\gamma  m, Ad \rho(\gamma) \cdot v)$
with $\gamma \in \pi_1(\M)$ acting on $\tilde \M$
by covering transformations and on ${\g}$
by the adjoint action of the holonomy $\rho(\gamma)$.
At each point $p$ of $\tilde \M$, the fiber of
$\tilde E$ splits as a direct
sum of infinitesimal pure translations and infinitesimal pure
rotations about $p$; these can be identified with
$T_p \tilde\M$ and $so(n)$ respectively.

We now lift $v$ to a section $s : \tilde \M \to \tilde E$
by choosing an ``osculating'' infinitesimal isometry $s(m)$ which best
approximates the vector field $v$ at each point $m \in \tilde \M$.
Thus $s(m)$ is the unique infinitesimal isometry
whose translational part and rotational part
at $m$ agree with the values of $v$ and $curl~v$ at $m$.
(This is the ``canonical lift'' as defined in \cite{HK1}.)
In particular, if $v$ is itself an infinitesimal isometry of
$\tilde \M$ then $s$ will be a constant function.

Using the equivariance property of the developing maps
it follows that $s$ satisfies an ``automorphic'' property:
$s(\gamma m) - Ad \rho(\gamma) s(m)$
is a {\em constant} infinitesimal isometry,
given by the variation $\dot \rho(\gamma)$ of holonomy isometries
$\rho_t(\gamma) \in G$ (see Prop 2.3(a) of \cite{HK1}).
Here $\dot \rho: \pi_1(\M) \to \g$
satisfies the cocyle condition $\dot \rho(\gamma_1 \gamma_2) = \dot
\rho(\gamma_1)+ Ad \rho (\gamma_1) \dot \rho(\gamma_2)$,
so represents a class in group cohomology
$[\dot \rho] \in H^1(\pi_1(\M);Ad \rho)$,
describing the variation of holonomy representations $\rho_t$.

Regarding  $s$ as a vector-valued function with values
in the vector space $\g$, its differential $\tilde\omega = ds$
satisfies $\tilde\omega(\gamma m) =  Ad \rho(\gamma) \tilde\omega(m)$ so
it descends to a closed 1-form $\omega$ on $\M$ with values in the bundle $E$.
Hence it determines a de Rham cohomology
class $[\omega]\in H^1(\M; E)$.
This agrees  with the cohomology class $[\dot \rho]$
under the de Rham isomorphism $H^1(\M; E) \cong H^1(\pi_1(\M);Ad \rho)$.
Also, we note that the translational part of $\omega$
can be regarded as
a $T\M$-valued $1$-form on $\M$. This is exactly the
form $\eta$ defined above (see Prop 2.3(b) of \cite{HK1}),
describing the infinitesimal change in metric on $\M$.

On the other hand, a family of hyperbolic structures determines
only an equivalence class of families of developing maps and
we need to see how replacing one family by an equivalent family
changes the cocycles.  Recall that a family equivalent to $\Phi_t$ is
of the form $k_t \Phi_t \tilde f_t$ where $k_t$ are isometries of $\H^n$ and
$\tilde f_t$ are lifts of diffeomorphisms $f_t$ from $\M$ to itself.
We assume that $k_0$ and $\tilde f_0$ are the identity.

The $k_t$ term changes the path $\rho_t$ of holonomy representations
by conjugating by $k_t$.  Infinitesimally, this changes the
cocycle $\dot \rho$ by a coboundary in the sense of group cohomology.
Thus it leaves the class in $H^1(\pi_1(\M);Ad \rho)$ unchanged.
The diffeomorphisms $f_t$ amount to choosing a different map
from $\M_0$ to $\M_t$. But $f_t$ is isotopic
to $f_0 = {\rm identity}$, so the lifts $\tilde f_t$
don't change the group cocycle at all.  It follows that equivalent families
of hyperbolic structures determine the same group cohomology class.

If, instead, we view the infinitesimal deformation as represented
by the $E$-valued $1$-form $\omega$, we note that the infinitesimal
effect of the isometries $k_t$ is to add a constant to
$s: \tilde \M \to \tilde E$.  Thus, $d s$, its
projection $\omega$, and the infinitesimal variation of metric are
all unchanged.  However, the infinitesimal
effect of the $\tilde f_t$ is to change the vector field on
$\tilde \M$ by the lift of a globally defined vector field
on $\M$.  This changes $\omega$ by the derivative of a
{\em globally defined} section of $E$.  Hence, it doesn't
affect the de Rham cohomology class in $H^1(\M; E)$.
The corresponding infinitesimal change of metric is altered
by the Lie derivative of a globally defined vector field
on $\M$.

Since, within an equivalence class of an infinitesimal deformation,
we are free to choose an identification of $\M_0$ with $\M_t$,
we can try to find a canonical choice with particularly nice
analytic properties.  A natural choice would be a harmonic map.
At the infinitesimal level, this corresponds to choosing a
Hodge representative for the de Rham cohomology class in $H^1(\M; E)$.
The translational part, which describes the infinitesimal change
in metric, is a {\em harmonic} $T\M$-valued 1-form.  These are studied
in detail for the case of cone-manifolds in \cite{HK1}.
They correspond to variations of metric which are
$L^2$-orthogonal to the trivial variations given by
the Lie derivative of compactly supported vector fields on $\M$.

One special feature of the 3-dimensional case is the
{\em complex structure} on the Lie algebra
$\g \cong sl_2 \C$ of infinitesimal isometries of $\H^3$.
The infinitesimal rotations fixing a point $p \in \H^3$
can be identified with $su(2) \cong so(3)$, then
the infinitesimal pure translations at $p$
correspond to $i \, su(2) \cong T_p \H^3$. Geometrically,
if $t \in T_p \H^3$ represents an infinitesimal translation,
then $i t $ represents an infinitesimal rotation with
axis in the direction of $t$. Thus,
on a hyperbolic 3-manifold $\M$ we can identify
the bundle $E$ of (germs of) infinitesimal isometries
with the {\em complexified} tangent bundle $T\M \otimes \C $.

We now specialize to the case of interest in this paper,
3-dimensional hyperbolic cone-manifolds.  We recall some of
the results and computations derived in \cite{HK1}.  The reader
is referred to that paper for further details.

Let $M_t$ be a smooth family of hyperbolic cone-manifold
structures on $M$ with cone angles $\alpha_t$ along $\Sigma$, where
$0 \le \alpha_t  \le 2\pi$. Note that, locally, $M_t$ is uniquely
determined by $\alpha_t$, by the local rigidity results of \cite{HK1}.
Let $U=U_R$ denote an embedded tube consisting of points distance at
most $R = R_t$ from the singular locus $\Sigma$.

By the Hodge theorem proved in \cite{HK1},
the infinitesimal deformation of hyperbolic structures
(``${d\over dt}(M_t)$'') can be represented
by a unique harmonic
$T\M$-valued 1-form $\eta$ on $\M = M-\Sigma$
such that 
$$D^* \eta =0,~ D^* D \eta = -\eta,$$
where $D$ is the exterior covariant derivative on such forms and
$D^*$ is its adjoint. In addition,
$\eta$ and $D\eta$ are symmetric and traceless, 
and inside $U$ we can write
$$\eta = \eta_0 + \eta_c$$
where $\eta_0$ is a ``standard'' (non-$L^2$) form,
and $\eta_c$ is a correction term with $\eta_c$,
$D\eta_c$ in $L^2$. Further, only $\eta_0$ changes
the holonomy of the meridian and longitude on the torus
$T_R = \bd U_R$. 

Alternatively, we can represent the
infinitesimal deformation by the 1-form with values
in the infinitesimal local isometries of $\M$:
\begin{eqnarray} 
\omega = \eta + i *\! D\eta.
\label{omega}\end{eqnarray}
 There is an analogous decomposition of $\omega$
in the neighborhood $U$ as $\omega = \omega_0 + \omega_c$
where only $\omega_0$ changes the holonomy and $\omega_c$ is in $L^2$.

The tubular neighborhood $U$ of the singular locus 
will be mapped by the developing
map into a neighborhood in ${\H}^3$ of a geodesic.  If we use 
cylindrical coordinates, $(r,\t, \z)$, the hyperbolic metric is
$dr^2 + \sinh^2  r\, d\t^2 + \cosh^2 r\, d\z^2$,
where the angle $\theta$ is defined modulo the cone angle $\alpha$.
We denote the  moving co-frame adapted to this coordinate system by
$(\omega_1, \omega_2, \omega_3) = (dr, \sinh r\, d\t, \cosh r\, d\z)$.

To define our standard forms, we use the cylindrical coordinates
on $U$ defined above, and we denote by
$e_1, e_2, e_3$ the orthonormal frame in $U$ 
dual to the co-frame $\omega_1,
\omega_2, \omega_3$.  
In particular, $e_2$ is tangent to the meridian and 
$e_3$ is tangent to the singular locus, which is homotopic in the 
cone-manifold to the longitude. 
We can interpret an $E$-valued 1-form as a complex-valued section of 
$T\M \otimes T^*\M \cong Hom(T\M,T\M)$. 
Then an element of
$T\M \otimes T^*\M$ can be described as a matrix whose 
$(i,j)$ entry is the coefficient of $e_i \otimes \omega_j$.

Explicitly, $\omega_0$ is a linear combination of
the forms given in (23) and (24) of \cite{HK1}.
The form $\omega_m = \eta_m + i *\! D \eta_m$ below 
is a ``standard'' closed and co-closed (non-$L^2$) form
which represents an infinitesimal deformation which
decreases the cone angle but does not change the real 
part of the complex length of the meridian.
It preserves the property that the meridian
is elliptic and, hence, that there is a cone-manifold structure.

\begin{eqnarray} \omega_m = \displaystyle{
 \begin{bmatrix}
 {-1\over\cosh^2(r) \sinh^2(r)}& 0& 0\\
 0& {1\over \sinh^2(r)}&{-i\over\cosh(r) \sinh(r)} \\
 0 & {-i\over\cosh(r) \sinh(r)} &{-1\over\cosh^2(r)} 
 \end{bmatrix} }\label{eq1}\end{eqnarray}

 The form $\omega_l = \eta_l + i *\! D\eta_l$
 below is a ``standard'' closed and co-closed, 
 $L^2$ form which stretches the 
 singular locus, but leaves the holonomy of the meridian 
 (hence the cone angle) unchanged. 
 
\begin{eqnarray}\omega_l = \displaystyle{
  \begin{bmatrix}
{-1\over\cosh^2(r)}& 0& 0\\
 0& -1&{-i \sinh(r)\over\cosh(r)} \\
 0 & {-i \sinh(r)\over\cosh(r)} &{\cosh(r)^2 + 1\over \cosh(r)^2}
 \end{bmatrix} } \label{eq2}\end{eqnarray}

The effect of $\omega_m$ and $\omega_l$ on the complex lengths of
the group elements on the boundary torus was computed in \cite{HK1}
(pages 32-33). For a detailed explanation for 
these computations we refer to this reference. 
We merely record the results here. 

\begin{lemma} \label{dlengths} The effects of the infinitesimal deformations given
by the standard forms on 
the complex length, ${\cal L}$, 
of any peripheral curve are as follows.
\begin{itemize} 
\item[(a)] For $\omega_m$,
\begin{eqnarray*} {d \over dt} ({\cal L})  = 
-2{\cal L}.\end{eqnarray*}
\item[(b)] For $\omega_l$,
\begin{eqnarray*} {d \over dt} ({\cal L})  = 
  2 \Re({\cal L}),\end{eqnarray*}
where $\Re({\cal L})$ denotes the real length of the curve. 
\end{itemize}
\end{lemma}

\begin{remark}
A {\em meridian} curve has complex length $i \alpha$. So the
effect of $\omega_m$ on its derivative is $ - 2 i\alpha$. 
This shows that the meridian 
remains elliptic and that the derivative of the cone angle $\alpha$ is 
$- 2 \alpha $. Similarly, for $\omega_l$, the complex length of the meridian 
has derivative zero.

If ${\cal L}$ denotes the complex length of the {\em longitude},  
then the real part of ${\cal L}$ is the length $\l$ of the singular locus. 
Thus for $\omega_m$, the derivative of $\l$ is $-2 \l$. For $\omega_l$,
the derivative of $\l$ is $2\l$.  

\end{remark}

The infinitesimal changes in the complex lengths of the elements 
of the fundamental group of the torus uniquely determine a complex 
linear combination of $\omega_m$ and $\omega_l$ and conversely
any such linear combination determines the infinitesimal
changes in these complex lengths.  The coefficient of $\omega_m$ uniquely 
determines and is determined by the change in the meridian since 
$\omega_l$ leaves the complex length of the meridian unchanged.
By our computations above the length of the meridian remains
pure imaginary (i.e. an elliptic element) precisely when the
coefficient is real.

The smooth family of structures $M_t$ is determined by a
choice of parametrization of the cone angles $\alpha_t$
and we are free to choose this as we wish.
The value of the coefficient for $\omega_m$ is determined by 
the derivative of  the cone angle.
It turns out to be useful to parameterize the cone-manifolds by the
{\em square} of the cone angle; i.e., we will let $t = \alpha^2$.
Since the  derivative of the square of the cone angle is $1$
and the derivative of $\alpha$ under $\omega_m$ is $-2 \alpha$,
we have
\begin{eqnarray}
\omega_0 ={-1\over 4\alpha^2} \omega_m + (x+iy) \omega_l
\label{omega_0}\end{eqnarray}
for some real constants $x$ and $y$.
One of the goals of this section is to estimate the values of
$x$ and $y$.  This will allow us to estimate the infinitesimal change
in all of the complex lengths of curves on the torus. In
particular, we can estimate the change in the length of the
singular locus.

The estimates in this section can
be viewed as effective versions of local rigidity arguments.
The basic idea behind local rigidity is to represent an
infinitesimal deformation by a harmonic representative in
the cohomology group $H^1(\M;E)$. 
The symmetric real part of this representative is a 1-form with values
in the tangent bundle of $\M$.  
Harmonicity, and the fact that it will be volume preserving 
(this takes a separate argument), imply that the 
1-form satisfies a Weitzenb\"ock-type formula:
$$D^*D\eta ~+~ \eta ~=~ 0$$
where $D$ is the exterior covariant derivative on such forms and
$D^*$ is its adjoint.  Taking the $L^2$ inner product of this
formula with $\eta$ and integrating by parts gives the formula
$$ ||D\eta||_X^2 ~+~ ||\eta||_X^2 ~=~ 0$$
when $\M$ is closed.  (Here $||\eta||_X^2$ denotes the square of the 
$L^2$ norm of $\eta$ on $X$. The pointwise $L^2$ norm is denoted simply
by $||\eta||$.)  Thus $\eta ~=~ 0$ and the deformation is
trivial. This is the proof of local rigidity for
closed hyperbolic 3-manifolds.

When $\M$ has boundary or is non-compact, there will be a boundary
term $\b$:
$$ ||D\eta||_X^2 ~+~ ||\eta||_X^2 ~=~ \b .$$ 
If the boundary term is non-positive, the same conclusion
holds: the deformation is trivial.  When $\M = M - \Sigma$, where
$M$ is a hyperbolic cone-manifold with cone angles at most $2\pi$ along its singular
set $\Sigma$, it was shown in \cite{HK1} that, for
a deformation which leaves the cone angle fixed, it is possible
to find a representative as above for which the boundary term
goes to zero on the boundary of tubes around the singular locus whose
radii go to zero.  Again, such an infinitesimal deformation must be trivial.
This proves local rigidity rel cone angles.

The argument for local rigidity rel cone angles actually shows that
the boundary term is negative when the cone angle is unchanged.  Note
that leaving the cone angle unchanged is equivalent to the vanishing of the
coefficient of $\omega_m$.  As we shall see below the boundary term
for $\omega_m$ by itself is positive.  Roughly speaking, $\omega_m$
contributes positive quantities to the boundary term, while everything
else gives negative contributions. (There are, of course, also some cross-terms.)
We think of ${-1\over 4\alpha^2} \omega_m $ as a preliminary model for
the infinitesimal deformation in a tube around the singular locus.  Then
this is ``corrected" by adding $(x+iy) \omega_l$ to get the actual change
in complex lengths and then by adding a further term $\omega_c$ 
that doesn't change the holonomy at all. 
The condition
that the boundary term for the actual representative (model plus
the other terms) be positive puts strong restrictions on
these ``correction" terms.  This is the underlying philosophy for 
the estimates in this section.

In order to implement these ideas we need to 
derive a formula for the boundary term.
For details we refer to \cite{HK1}.

The Hodge Theorem (\cite{HK1}) for cone-manifolds gives a closed and
co-closed $E$-valued form $\omega = \eta + i *\!D\eta$ satisfying
$D^{*}D \eta = - \eta$. 
Integration by parts, as in \cite{HK1} (Proposition 1.3 and p. 36), 
over any sub-manifold $N$ of $\M$ with boundary $\bd N$ gives:

\begin{lemma} For any closed and co-closed form 
$\omega = \eta + i*\! D\eta$ satisfying $D^{*}D \eta = - \eta$, 
and any submanifold $N$ with boundary $\bd N$ oriented by the
{\em outward} normal,
\begin{eqnarray} 0 = \int_N ( ||\eta||^2 + ||\!*\!D\eta||^2) +  
\int_{\bd N} * D\eta \wedge \eta. \label{bdryeq}\end{eqnarray}
\end{lemma}

Note that in these integrals, $\alpha \wedge \beta$ denotes 
the real valued 2-form obtained
using the wedge product of the form parts, and the geometrically
defined inner product on vector-valued parts.

Denote by $U_r$ the tubular neighborhood of points at distance
less than or equal to $r$ from the singular locus. It will always
be assumed that $r$ is small enough so that $U_r$ will be embedded. 
Let $T_r$ denote the boundary torus of $U_r$, 
oriented with ${\bd\over \bd r}$ as outward normal. 
We define
\begin{eqnarray}\b_r(\alpha,\beta) = \int_{T_r} *D\alpha \wedge \beta.
\label{b_defn}\end{eqnarray}
We emphasize that $T_r$ is oriented as above, so that
$\omega_2 \wedge \omega_3 = \sinh r \, \cosh r \, d\theta \wedge d\z$
is the oriented area form.

Fix a value $R$ for the radius and let $N = X - U_R$. 
Then $\bd N = -T_R$, where the minus sign denotes the opposite
orientation (since $-{\bd\over \bd r}$ is the outward normal
for $N$).  Applying (\ref{bdryeq}) in this case, we obtain

\begin{corollary} Let $N = X - U_R$ be the complement
of the tubular neighborhood of radius $R$ around the singular locus.  
Then, for any closed and co-closed form $\omega = \eta + i*\!D\eta$ 
satisfying $D^{*}D \eta = - \eta$,
\begin{eqnarray}  \b_R(\eta,\eta) =  ||\eta||_N^2 + ||\!*\!D\eta||_N^2 
= ||\omega||_N^2 .
\label{etabdry}\end{eqnarray}
\end{corollary}

In particular, we see that the boundary term 
$b_R(\eta,\eta)$ is {\em non-negative}.
Writing $\eta = \eta_0 + \eta_c$ as before, 
we analyze the contribution from each part.
First, we note that the cross-terms vanish so that
the boundary term is simply the sum of two boundary terms:

\begin{lemma} \label{bdrydecomp} $\b_R(\eta,\eta) = \b_R(\eta_0,\eta_0) 
+ \b_R(\eta_c,\eta_c).$

\end{lemma}

\begin{proof} Expanding, we have that
$\b_R(\eta,\eta) = \b_R(\eta_0 + \eta_c,\eta_0 + \eta_c)
= \b_R(\eta_0,\eta_0) + \b_R(\eta_c,\eta_c) 
+ \b_r(\eta_0,\eta_c) + \b_r(\eta_c,\eta_0).$
So it suffices to show that $\b_r(\eta_0,\eta_c) = \b_r(\eta_c,\eta_0) = 0$.

This follows from the Fourier decomposition
for $\eta_c$ obtained in \cite{HK1}. The term 
$\eta_c$ is the infinitesimal change of metric induced by a 
vector field that satisfies a harmonicity 
condition in a neighborhood of the singular locus.  The main point
is that $\eta_c$ has no purely radial terms.  This can be seen 
from Proposition 3.2 of that paper, where the purely radial 
solutions correspond, in the notation used there, to the 
case $a=b=0$.  There is a $3$-dimensional solution space allowed 
by the chosen domain for the harmonicity equations 
(equations (21) in that paper). It becomes $2$-dimensional after 
the conclusion that the deformation is volume-preserving.
However, there is an obvious $2$-dimensional space of radial solutions
coming from the infinitesimal rotations and translations along the
axis corresponding to the singular locus.  Since these are isometries,
they don't contribute anything to the change of metric, $\eta_c$.

On the other hand, $\eta_0$ only depends on $r$ by definition,
so each term in the integrands for $\b_r(\eta_0,\eta_c)$ and
$\b_r(\eta_c,\eta_0)$ has a trigonometric factor which integrates
to zero over the torus $T_r$.
\end{proof}

Next, we show that the contribution, $\b_{R}(\eta_c,\eta_c)$,  
from the part of the ``correction term" that doesn't affect 
the holonomy is {\em non-positive}.  In fact,

\begin{lemma}
\begin{eqnarray} \b_{R}(\eta_c,\eta_c) = 
- ( ||\eta_c||_{U_R}^2 + ||\!*\!D\eta_c||_{U_R}^2) =
- ||\omega_c||_{U_R}^2.
\label{etacbdry}
\end{eqnarray}
\end{lemma}

\begin{proof} 

Consider a region $N=U_{r_1,r_2}$ in $U_R$
bounded by the tori $T_{r_1}$ and $T_{r_2}$ where  
$0< r_1 < r_2 \le R$.
Then $\bd N = T_{r_2} \cup -T_{r_1}$ where, as before,
the minus sign denotes the opposite orientation.

The equation (\ref{bdryeq}), applied to this region with 
$\eta = \eta_c$, gives 
$$0 =  \int_{U_{r_1,r_2}} ( ||\eta_c||^2 + ||\!*\!D\eta_c||^2) +
\int_{T_{r_2}} *D\eta_c \wedge \eta_c - \int_{T_{r_1}} *D\eta_c \wedge \eta_c,$$
or
\begin{eqnarray}
\b_{r_2}(\eta_c,\eta_c) - \b_{r_1}(\eta_c,\eta_c)
~=~ -\int_{U_{r_1,r_2}} ( ||\eta_c||^2 + ||\!*\!D\eta_c||^2). 
\label{r1r2bdry}
\end{eqnarray}

The main point here is that
$\lim_{r\to 0} \b_r(\eta_c,\eta_c) = 0$.
This is a restatement of the main result in section 3 of \cite{HK1},
since $\eta_c$ represents an infinitesimal deformation which
doesn't change the cone angle. 

Applying (\ref{r1r2bdry}), with $r_2 = R$ and taking the
limit as $r_1 \to 0$ we obtain the desired result.
\end{proof}

Combining Lemma \ref{bdrydecomp} with (\ref{etabdry}) and
(\ref{etacbdry}), we obtain:

\begin{eqnarray}
\b_R(\eta_0,\eta_0) = ||\omega||_{X-U_R}^2 + ||\omega_c||_{U_R}^2. 
\label{eta0bdry}\end{eqnarray}

In particular, this shows that 
\begin{eqnarray}\b_R(\eta_0,\eta_0) \ge 0, \label{eta0pos}\end{eqnarray} 

\smallskip
\begin{remark*}
This positivity is the only application of formula (\ref{eta0bdry}) we will
use in this paper.  However, we note here for future reference that
an upper bound on $\b_R(\eta_0,\eta_0)$ provides an upper bound on
the $L^2$ norm of $\omega$ on the complement of the tubular neighborhood
of the singular locus.  Such a bound can be used to bound the infinitesimal
change in geometric quantities, like lengths of geodesics, away from
the singular locus.  Similarly, an upper bound on $\b_R(\eta_0,\eta_0)$ 
provides an upper bound on the $L^2$ norm of the correction term $\omega_c$ 
in the tubular neighborhood itself.  This can be used to bound changes
in the geometry of the tubular neighborhood that are not detected simply
by the holonomy of group elements on the boundary torus.
\end{remark*}
%\smallskip

In the remainder of this section we will use the inequality (\ref{eta0pos})
to find bounds on the infinitesimal variation of the holonomy of the
peripheral elements.  Of particular interest will be bounding the
variation in the length of the singular locus (which equals the real part
of the complex length of any longitude of the boundary torus).
To this end, we further decompose $\eta_0$ as a sum of a component
that changes the cone angle and ones that leave it unchanged.

Recall that $\omega_0 ={-1\over 4\alpha^2} \omega_m + ( x + i y) \omega_l$
so 
$$\eta_0 = \Re( \omega_0) 
= {-1\over 4\alpha^2} \eta_m + x \eta_l - y *\!D\eta_l.$$  

The basic principle here is that the contribution of the 
$\eta_m$ term to $\b_R(\eta_0,\eta_0)$ is positive, while those
of the $\eta_l$ and $*D\eta_l$ terms are negative. (The
cross-terms only complicate matters slightly.)  The coefficient
of the $\eta_m$ term is fixed by the choice of parametrization
of the family of cone-manifolds by $t = \alpha^2$ so the fact
that $\b_R(\eta_0,\eta_0)$ is positive will provide a bound on the
coefficients $x$ and $y$.

We calculate
\begin{eqnarray*}
\b_R(\eta_0,\eta_0) &=& 
{1\over 16\alpha^4} \b_R( \eta_m, \eta_m)
+ x^2 \b_R( \eta_l, \eta_l) + y^2 \b_R(*D\eta_l,*D\eta_l) \\
&&-{x\over 4\alpha^2} ( \b_R(\eta_m,\eta_l) + \b_R(\eta_l,\eta_m))  
+{y\over 4\alpha^2} ( \b_R(\eta_m,*D\eta_l) + \b_R(*D\eta_l,\eta_m)) \\
&&-xy ( \b_R(\eta_l,*D\eta_l) + \b_R(*D\eta_l,\eta_l)) .
\end{eqnarray*}

Now, using the explicit formulas for $\eta_m$ and $\eta_l$, we find
\begin{eqnarray}\b_R( \eta_m, \eta_m) 
={1\over \sinh(R) \cosh(R)} 
	\biggl({1\over \sinh^2(R) } + {1\over \cosh^2(R)}\biggr) \area(T_R) ,
\label{bmm}\end{eqnarray}
\begin{eqnarray}\b_R( \eta_l, \eta_l) 
= \b_R(*D\eta_l,*D\eta_l)
= {-\sinh(R)\over \cosh(R)} 
 \biggl(  2 + {1\over \cosh^2(R)}\biggr) \area(T_R),
\label{bll}\end{eqnarray}
\begin{eqnarray}\b_R( \eta_m, \eta_l) 
= {-1\over \sinh(R) \cosh(R)} 
 \biggl(  2 + {1\over \cosh^2(R)}\biggr) \area(T_R),
\label{bml}\end{eqnarray}
\begin{eqnarray}\b_R( \eta_l, \eta_m) 
= {\sinh(R)\over \cosh(R)} 
	\biggl({1\over \sinh^2(R) } + {1\over \cosh^2(R)}\biggr) \area(T_R),
\label{blm}\end{eqnarray}
and the other terms vanish.

It simplifies matters slightly and is somewhat illuminating 
to rewrite the value of the boundary term $\b_R( \eta_0, \eta_0)$ 
using the geodesic length $m$ of meridian on the flat boundary torus $T_R$.
Recall that
$$m ~=~ \alpha \,\sinh (R).$$
Then we obtain
$$\b_R(\eta_0,\eta_0)/\area(T_R) = a (x^2+y^2) + b x + c ,$$
where
$$a =  {-\sinh(R)\over \cosh(R)} \biggl(  2 + {1\over \cosh^2(R)}\biggr) ~=~
-\tanh (R) ~ {2 \cosh^2(R) + 1 \over \cosh^2(R)},$$
$$b = {1\over 4\alpha^2} \biggl( {2\over \cosh^3(R) \sinh(R)}\biggr) ~=~ 
{1 \over m^2} ~~{\tanh(R) \over 2 \cosh^2(R)},$$
$$c = {1\over 16\alpha^4} {1\over \sinh(R) \cosh(R)} 
\biggl({1\over \sinh^2(R) } + {1\over \cosh^2(R)}\biggr) ~=~
{1 \over m^4} ~~{\tanh(R) + \tanh^3(R) \over 16}.$$

Completing the squares gives
\begin{eqnarray*}
\b_R(\eta_0,\eta_0)/\area(T_R) &=&  a (x^2 + y^2) + b x + c  \\
&=& a \biggl( \biggl(x+ {b\over 2a}\biggr)^2 + y^2 \biggr) 
+ { 4 a c - b^2 \over 4 a}.
\label{compsq}
\end{eqnarray*} 

By direct computation we see that
\begin{eqnarray}
4 a c - b^2  ~=~ {\tanh^2(R) \over m^4}.
\label{discrim}
\end{eqnarray}

Since $a$ is negative, we obtain the following estimate for the
boundary term $\b_R(\eta_0,\eta_0)$.  As noted before, we won't use
this estimate in this paper, but rather record it for future reference.
\begin{eqnarray}
\b_R(\eta_0,\eta_0)/\area(T_R) \le { 4 a c - b^2 \over 4 a}
~=~  {1\over 4 m^4} ~ {\sinh(R) \cosh(R) \over 2 \cosh^2(R) +1}.
\label{b00}\end{eqnarray}

Our main application of the positivity (\ref{eta0pos}) of 
$\b_R(\eta_0,\eta_0)$ is that, using (\ref{compsq}),
we can conclude that: 
\begin{eqnarray} \biggl(x+ {b\over 2a}\biggr)^2 + y^2  \le { b^2 - 4ac \over 4 a^2} 
={1 \over 4 m^4} ~ {\cosh^4(R) \over (2 \cosh^2(R) +1)^2}.
\label{etalineq}\end{eqnarray}

This implies, in particular, that $x$ lies in the interval of radius 
$${1 \over 2 m^2} ~ {\cosh^2(R) \over 2 \cosh^2 (R) +1}$$
around 
$${-b \over 2 a} ~=~ {1 \over 4 m^2} ~ {1 \over 2 \cosh(R)^2 +1}.$$
In other words,
$x$ lies in the interval $[x_1,x_2]$ where
$$x_1 = {-b - \sqrt{b^2-4 a c} \over 2a} = 
{-1 \over 4 m^2} ~ {2 \cosh^2(R) -1 \over 2 \cosh(R)^2 +1},$$
and 
$$x_2 = {-b + \sqrt{b^2-4 a c} \over 2a}  
= {1 \over 4 \m^2}.$$                                                                          

\begin{remark*} It is useful to rewrite the factor in
the formula above for $x_1$ as
$$ {2 \cosh^2(R) -1 \over 2 \cosh(R)^2 +1} ~=~ 
{ 2\sinh^2(R) + 1 \over 2 \sinh^2(R) + 3}.$$
Note that it is monotone increasing in $R$, 
taking on values between ${1\over 3}$ and $1$.
\end{remark*}

By Lemma \ref{dlengths} the effect of $\omega_0 $ on the complex length, ${\cal L}$, 
of any peripheral curve is given by
\begin{eqnarray} {d \over dt} ({\cal L})  = 
{-1\over 4\alpha^2} (-2{\cal L}) + (x+iy) (2 \Re({\cal L})),
\label{dLdt}\end{eqnarray}
where $\Re({\cal L})$ denotes the real length of the curve. 

In particular,  the derivative of the real length
$\l$ of the longitude 
(the length of the singular locus)
 satisfies
\begin{eqnarray}{d \l \over dt} = {\l \over 2 \alpha^2} ( 1 + 4\alpha^2 x).
\label{dldt}\end{eqnarray}
Since $t = \alpha ^2$, we conclude that
\begin{eqnarray} {d \l \over d \alpha} = {\l \over \alpha} ( 1 + 4\alpha^2 x).
\label{dldalpha}\end{eqnarray}

Putting this formula for the derivative of the length of the 
singular locus together with the estimates above for the 
coefficient $x$ (and recalling that $m = \alpha \sinh (R)$), 
we obtain the main result of this section:

\begin{theorem} \label{keyboundtheorem}
Consider any smooth family of cone structures on $M$, all of whose cone
angles are at most $2 \pi$.  For any component of the singular set, let
$\l$ denote its length and $\alpha$ its cone angle.  Suppose there is
an embedded tube of radius $R$ around that component.  Then
$${d \l \over d \alpha} = {\l \over \alpha} ( 1 + 4\alpha^2 x),$$
where 
\begin{eqnarray} {-1 \over \sinh^2(R)} 
	\biggl({ 2\sinh^2(R) + 1 \over 2 \sinh^2(R) + 3}\biggr)
 \le 4 \alpha^2 x \le {1\over \sinh^2(R) }.
\label{keybound}\end{eqnarray}
\end{theorem}

\smallskip

\begin{remark} \label{lincreases}
This implies that $\l$ is an increasing function of $\alpha$
provided the tube radius $R$ is large enough. Explicitly
$${d \l \over d\alpha} \ge 0$$ provided
$${1 \over \sinh^2(R)} 
\biggl({ 2\sinh^2(R) + 1 \over 2 \sinh^2(R) + 3}\biggr)\le 1$$
which simplifies to
$$R \ge \arcsinh({1\over \sqrt 2}) \approx 0.65848.$$
This has implications concerning the variation of the volume $V$ of
a family of cone-manifolds due to the Schl\"afli formula 
(see \cite{Ho}, \cite[Theorem 3.20]{CHK1}):
\begin{eqnarray} {d V \over d\alpha} = -{1\over 2} \l.
\label{schlafli}\end{eqnarray}
Since
$$ {d^2 V \over d\alpha^2 } = -{1\over 2}{d \l \over d\alpha} \le 0$$
for these values of $R$, the volume function will be a concave function of $\alpha$ as
long as the tube radius is sufficiently large.

More specifically, if one considers a family of cone-manifolds with a single 
component of the singular locus in which the cone angle is {\em decreasing}, the
total change, $\Delta V$, in the volume will be positive.  If
$R \ge \arcsinh({1\over \sqrt 2})$ throughout the deformation, then
we obtain the inequality
\begin{eqnarray} \Delta V \leq  {|\Delta \alpha| \over 2} \l_0,
\label{easyvolume}\end{eqnarray}
where $\Delta \alpha$ denotes the total change in cone angle and
$\l_0$ denotes the initial length of the singular locus.

\end{remark}

In Section \ref{tuberad}, we will see how to control the tube radius 
by controlling the length of the singular locus. This will lead to 
sharper estimates for the change in volume by integrating the more 
detailed estimates for ${d \l \over d \alpha}$ which are derived 
there.  However, it seems worthwhile to note that the above estimates
follow immediately from (\ref{keybound}).

\section{Geometric limits of cone-manifolds}\label{geomlim}

This section is primarily devoted to the proof of Theorem \ref{thm2}.

In general, the limiting behavior of a sequence of hyperbolic cone-manifolds
can be quite complicated. In particular, it can collapse to a lower
dimensional object or the singular locus can converge to something of
higher complexity.  However, by the results of Section \ref{tuberad},
we will be able to assume that there is a lower bound to the tube 
radius around each component of $\Sigma$ and that the geometry of the
boundary of the tube doesn't degenerate.  This greatly simplifies
matters, essentially reducing us to the manifold case.

Given a sequence of hyperbolic cone-manifold structures $M_i$ on $(M,\Sigma)$,
remove disjoint, embedded equidistant tubes around each component
of $\Sigma$.  The result is a sequence of smooth, hyperbolic
manifolds $N_i$ with torus boundary components, each of which
has an intrinsic flat metric.  Furthermore, the principal normal curvatures
are constant on each component, equalling $\kappa, {1 \over \kappa}$
(we assume that $\kappa \geq 1$).  When $\kappa > 1$ the lines of curvature are
geodesics in the flat metric corresponding to the meridional and
longitudinal directions, respectively. Note that the normal curvatures
and the tube radius, $R$, are related by ${\rm coth}~R = \kappa$ so
they determine each other.

We now formalize the structure of this type of boundary torus.  Let
$\H^3_R$ denote $3$-dimensional hyperbolic space minus the open
tube of points distance less than $R$ from a geodesic.  We allow the
values $0<R\leq \infty$, where $\H^3_{\infty}$ denotes the complement
of an open horoball based at a point at infinity.  We say that
a torus boundary component of a hyperbolic $3$-manifold is
{\em locally modelled on} $\H^3_R$ if, for some fixed $R$,
each point on the boundary torus has a neighborhood isometric to
a neighborhood of a point on the boundary of $\H^3_R$.  The
overlap maps are required to be restrictions of $3$-dimensional 
hyperbolic isometries.
This is equivalent to the condition that the torus have an induced
flat metric and have normal curvatures and lines of curvature as
in the previous paragraph.  Note that normal curvatures all equal to
$1$ corresponds to the case $R=\infty$.

\begin{definition} A hyperbolic $3$-manifold is
said to have {\em tubular boundary} if its
boundary consists of tori that are each locally modelled on
$\H^3_R$ for some $0<R\leq \infty$.  (The value of $R$ is
allowed to be different on different components of the boundary.)
\end{definition}

As noted above, one way these arise is when one removes tubular neighborhoods 
of the components of the singular locus of a cone-manifold.  On the other hand, 
we will see below that there is a canonical way to fill in any tubular 
boundary component. If a hyperbolic $3$-manifold with tubular boundary 
came from a hyperbolic cone-manifold by removing tubular neighborhoods, 
the filling process recovers the same cone-manifold.

To see this, first note that when $R = \infty$ the boundary torus has
all normal curvatures equal to $1$, so it can be identified with
a horosphere modulo a group of parabolic isometries fixing the corresponding
point at infinity.  This group action extends canonically to an action
on the horoball bounded by the horosphere.  In this case, the boundary is 
``filled in" with a cusp.  This is interpreted as a cone-manifold structure
with cone angle $0$.  If the tubular boundary came from removing a tubular
neighborhood of the ``singular locus", it must actually have been a cusp 
because the normal curvatures all equal $1$.  Furthermore, since the 
structure of the cusp is determined by the flat structure on the
boundary, the cusp replaced must be isometric to the one removed.  

To analyze the case of finite $R$, we note that the universal cover of 
the complement of a geodesic in $\H^3$ is isometric to $\R^3$ with metric in 
cylindrical co-ordinates $(r,\theta,\zeta)$, where $0 < r$,  given by
\begin{eqnarray} dr^2 ~+~ \sinh^2 r \, d\theta^2 ~+~ \cosh^2 r \, d\zeta^2. 
\label{cylmetric} \end{eqnarray}

A neighborhood of the tubular boundary is given by dividing out a
neighborhood of the plane $r=R$ in $\R^3$ by a $\Z \oplus \Z$ 
lattice in the $(\theta, \zeta)$-plane.  The above metric descends to the metric
in a neighborhood of the tubular boundary.  In particular, the boundary
is the image of $r=R$ and the principal curvatures, $\kappa, 1/\kappa$, 
are in the $\theta$, $\zeta$ directions, respectively.  
The metric on the tubular boundary can be canonically
extended by adding the quotient of the region $r \in (0,R]$ by the
$(\theta, \zeta)$ lattice group.  This metric is incomplete.  In general
its completion is singular, resulting in a hyperbolic structure 
``with Dehn surgery singularities" (see Thurston \cite{thnotes} 
for further discussion).  This structure includes cone-manifolds as a special case.
We will not be concerned with the more general type of singularity here, but
rather see below that the cone-manifold structures can be identified from
the structure on the tubular boundary.

If one removes a tubular neighborhood of a component of the singular
locus of a cone-manifold with cone angle $\alpha$, the boundary torus has a closed
geodesic in the meridian ($\zeta = {\rm constant}$) direction which is the
boundary of a totally geodesic, singular disc with cone angle $\alpha$
perpendicular to the core geodesic.  Conversely, we claim that if there is such a
closed geodesic, the completion defined above will be a cone-manifold.  To see this,
note that there is a closed meridian on the boundary torus
if and only if the lattice in $(\theta, \zeta)$ can be chosen to have one
generator of the form $(\alpha, 0)$.  We denote by $(\tau, \l)$ the other generator,
where necessarily $\l \neq 0$.  This corresponds to the first generator
being a rotation by angle $\alpha$ around the removed geodesic.  The
second generator translates distance $\l$ along the removed geodesic 
and rotates by angle $\tau$; i.e., it has complex length $\l + i \tau$.
Then the completion is obtained by adding in the quotient of the removed geodesic
(corresponding to $r=0$) under the action.

This is easily seen to be a cone-manifold
with cone angle $\alpha$.  In particular, the singular locus is the geodesic
added in the completion and the meridian, which is a closed geodesic in the 
flat metric on the original tubular boundary, bounds a singular, totally 
geodesic disk intersecting the singular locus in a single point.  
The flat structure on the tubular boundary can be constructed by taking
a flat cylinder of circumference $\m$ and height $\h$ and attaching it with
a twist of distance $\tw$.  The cone angle, $\alpha$, and the complex length 
$\l ~+~ i\tau$ are related to these quantities by the equations:
\begin{eqnarray*} \m &=& \alpha \sinh R, \\
\h &=& \l  \cosh R,\\
\tw &=& \tau  \sinh R.
\end{eqnarray*}

This implies that the region added is canonically determined
by the geometry of the boundary torus, the value of $R$, and
the fact that there is a closed geodesic in the meridian (principal 
curvature $\kappa > 1$) direction.  Thus, if the tubular boundary
structure arose from removing a tubular neighborhood of a component 
of the singular locus of a cone-manifold, the filling in process 
would recover the same cone-manifold structure.

The results proved in this section concern bilipschitz limits of
sequences of hyperbolic manifolds with tubular boundary. The above analysis
implies that if the members of the sequence all arise from cone-manifold
structures, and if the limit is a hyperbolic manifold
with tubular boundary, then it can be filled in to be a cone-manifold
also, and the results can be viewed in terms of bilipschitz
limits of cone-manifolds.

There are two advantages to considering sequences of hyperbolic
structures with such boundary data rather than studying sequences of
hyperbolic cone-manifolds directly.  First, the analysis of geometric 
limits is much simpler in the manifold setting.  Though the boundary
does introduce complications similar to those that arise for cone-manifolds,
it is easier to isolate them if singular locus is removed.
Secondly, the results of this section will apply to more general
singular structures than cone-manifolds.  In particular, they will
apply to a sequence of hyperbolic structures with Dehn surgery
singularities as long as there is a lower bound to the radii of disjoint
tubes around the singularities.  We expect to use this application
in a future paper.

A topological ball in a hyperbolic manifold with tubular boundary
will be called {\em standard} if it is isometric to a ball of radius 
$r >0$ in $\H^3$ or to a ball of radius $r > 0$ about a point on the 
boundary of $\H^3_R$.  In the latter case, we further require that 
$r<R$. This corresponds to the geometric condition that if the 
tube of radius $R$ were added back to $\H^3_R$ and the ball 
extended to to a ball in $\H^3$, then the extended ball would 
be disjoint from the geodesic core of the added tube.

The {\em injectivity radius} at a point $x$ in a
hyperbolic manifold, $N$, with tubular boundary is
$$inj(x,N)= \sup \{ r \mid B_r(x) \subset {\rm a~standard
~ball~in~}N\}.$$
Here $B_r(x)$ simply denotes the set of points in $N$ distance 
less than $r$ from $x$; there is no assumption on its topology.
We will write $inj(N)$ to denote $\inf_{x \in N}(inj(x,N))$.

Note that we {\em do not} assume that the standard
neighborhood is {\em centered} at the point $x$. This is
to avoid difficulties near the boundary: a point $x$
near, but not on, the boundary has only a small standard ball
centered at $x$, with radius at most the distance to the boundary.
However, there may be much larger standard balls which contain $x$ that
are centered at a point on the boundary.

It is important also to notice that because of the condition that
$R>r$ for a standard ball of radius $r$ centered at a point on a
boundary torus locally modelled on $\H^3_R$, a lower bound on
the injectivity radius of $N$ implies a lower bound on the
tube radii of all the boundary components.

The goal of this section is to find conditions on a family of hyperbolic
$3$-manifolds with tubular boundary that ensure that they converge
to a diffeomorphic manifold with such a structure.
The notion of convergence that we will use is based on
a distance between metric spaces defined using bilipschitz mappings.

\begin{definition} The {\em bilipschitz distance} between two metric
spaces $X,Y$ is the infimum of the numbers
\begin{eqnarray}
|{\rm log}\, {\rm lip}(f)| ~+~ |{\rm log}\, {\rm lip}(f^{-1})|
\end{eqnarray}
where $f$ ranges over all bilipschitz mappings from $X$ to $Y$
and ${\rm lip}(f)$ denotes the lipschitz constant of $f$.
\end{definition}

The bilipschitz distance between $X$ and $Y$ is defined to be $\infty$
if there is no bilipschitz map between them.  In particular, metric
spaces that are a finite distance apart are necessarily homeomorphic.
It is not hard to show that two compact metric spaces are bilipschitz
distance $0$ apart if and only if they are isometric.

For non-compact spaces, bilipschitz distance is not very useful because
it is so often infinite.  For many purposes, it is important to allow
a more flexible idea of convergence of sequences of metric spaces
than that induced simply by bilipschitz distance.
To make this idea precise, it is necessary
to choose a basepoint in each metric space.

\begin{definition} A sequence, $\{(Y_i, y_i)\}$, of metric spaces
with basepoint {\em converges} to $(Y,y)$ in the
{\em pointed bilipschitz topology} if, for each fixed $R>0$, the
radius $R$ neighborhood of $y_i$ in $Y_i$
converges with respect to the bilipschitz distance to the
radius $R$ neighborhood of $y\in Y$.
\end{definition}
Note that with this notion of convergence, a sequence of compact spaces
can converge to a non-compact space.  In particular, there is no
requirement that the $Y_i$ in a convergent sequence be eventually homeomorphic.
Convergence in the pointed bilipschitz topology means that the metric
spaces are becoming closer and closer to being isometric on
larger and larger diameter subsets.  However, when there is a uniform bound
to the diameter of all the spaces in the sequence, convergence is
independent of the choice of basepoint and is just convergence with
respect to the bilipschitz metric.

Our beginning point in the study of convergence of hyperbolic $3$-manifolds
with tubular boundary is a seminal and general theorem due to
Gromov.  It says that, under very mild conditions, (pinched curvature
and bounded injectivity radius at the basepoint), a sequence of complete,
pointed Riemannian manifolds will have a convergent subsequence in
this topology.  This theorem is actually a corollary of an even broader
compactness theorem, involving a much more general notion of convergence of
metric spaces, usually referred to as Gromov-Hausdorff convergence.
However, Gromov shows that, when applied to various classes of Riemannian
manifolds, this general notion of convergence implies convergence
in the pointed bilipschitz topology.  We will not need to use
the concept of Gromov-Hausdorff convergence in this paper, but rather
begin with its application to Riemannian manifolds.

\begin{theorem} \cite[Theorem 8.25]{GrFrench}, \cite [Theorem 8.20] {GrEnglish}
\label{Glipconv}
Consider a sequence of complete, pointed Riemannian manifolds $(N_i, v_i)$ with
pinched sectional curvatures $|k|\le K$ and injectivity radius at the basepoint, 
$v_i$, bounded below by $c > 0$.  Then there is a pointed Riemannian
manifold $(N,v)$ and a
subsequence of the $(N_i,v_i)$ which converges in the pointed bilipschitz topology
to $(N,v)$.  Furthermore, if there is a $D>0$ so that the diameters of the
$N_i$ are less than $D$ for all $i$, then the $N_i$ in the convergent 
subsequence will be diffeomorphic to $N$ for $i$ sufficiently large.
\end{theorem}

The fact that convergence in the metric is only lipschitz means that, a priori,
the limit metric is only $C^{0}$.  In \cite{GrFrench} and \cite{GrEnglish},
it is explained how a somewhat higher level of regularity can be achieved by
considering harmonic co-ordinates.  For closed manifolds, a complete proof along the 
lines sketched there appears in \cite{katsuda}.  Proofs along somewhat different lines appear
in \cite{Greene-Wu} and \cite{Peters}; these references also provide simple
examples showing why the limit metric won't be $C^2$ in general.  However,
if all the metrics in the sequence are of a special type, much stronger conclusions
are possible.  As explained in \cite{petersen}, p.307, if the approximating metrics
are Einstein, then use of the Einstein equation and elliptic regularity allows one
to bootstrap the regularity of convergence to any number of derivatives and the limit
metric will also be Einstein.

In our situation with constant curvature, things are vastly simpler.  The regularity
issues discussed above are all local.  The regularity of the convergence and of the
limit metric follow from local analysis on embedded balls of fixed radius.  In general, 
simply bounding the injectivity radius and curvature of a sequence of metrics
does not bound derivatives of the curvature and smoothness may be lost in the limit,
even locally.  However, since all metric balls of a fixed radius in hyperbolic $n$-space
are isometric, the bilipschitz limit of a sequence of hyperbolic $n$-balls of 
fixed radius will automatically be hyperbolic.  Thus, in the theorem above,
if the approximating manifolds are all hyperbolic, the limit manifold will be also.

The fact that we are considering manifolds with boundary means that we can't immediately
apply Theorem \ref{Glipconv} above.  Indeed, a few extra conditions on the boundary
are necessary, for example, to keep the boundary from collapsing to a point or to keep
two components on the boundary from colliding in the limit.  This has been worked out
in \cite{Kodani}, where Gromov's theorem is extended to manifolds with boundary
if one has the added conditions that the principal curvatures and intrinsic diameters 
of the components of the boundary are bounded above and below and that there is a 
lower bound to the width of an embedded tubular neighborhood of the boundary.  
We see in the proof below that, with our definition of injectivity radius, these 
conditions hold for manifolds with tubular boundary if the injectivity radius is 
bounded below for points on the boundary and the volume of the entire manifold 
is bounded above.

\begin{theorem} \label{tubelipconv} Let $(N_i, v_i)$ be a sequence of
hyperbolic $3$-manifolds with tubular boundary with basepoints $v_i$ on $\bd N_i$.
Assume there are constants $c, V >0$ such that, for all $i$,
$inj(x, N_i) \ge c$ for all $x \in \bd N_i$ and $\vol(N_i)<V$.
Then there is subsequence converging in the pointed bilipschitz topology
to a pointed hyperbolic $3$-manifold with tubular boundary,
$(N_\infty, v_{\infty})$.
Furthermore, if the diameters of all the $N_i$ are uniformly bounded, then
all the $N_i$ in the subsequence will be diffeomorphic to $N_\infty$ for
sufficiently large $i$.
\end{theorem}

\begin{remark} The bound on the volume will only be used to conclude that the
intrinsic diameters of the boundary components of all the $N_i$'s are uniformly
bounded.  Thus, the theorem remains true with the volume condition replaced by
such a bound on these intrinsic diameters.
\end{remark}

\begin{proof}
In order to apply the generalization in \cite{Kodani} of Theorem \ref{Glipconv}
we need to check the required conditions on the boundary.  Recall that points 
on the tubular boundary are locally modelled on $\H^3_R$ and that the definition
of injectivity radius implies that $inj(x, N_i)< R$ for such points.  Since
the principal curvatures on the boundary equal $\kappa, 1/\kappa$, where
$\kappa = \coth R$, a lower bound on the injectivity radius for boundary points
immediately bounds the principal curvatures above and below.

The definition of injectivity radius at a point requires that there will be a 
standard ball containing the set of points distance $r$ from the point, for 
any $r$ less than the injectivity radius.  The radius of the 
standard ball must be at least equal to this $r$.  But any standard ball 
in $\H^3_R$ containing a boundary point must be centered at some (possibly
different) point on the boundary of $\H^3_R$. This implies that there
is a tubular neighborhood around the boundary with a lower bound to its width.

Finally, we need to see that the intrinsic diameters of the boundary components
are bounded above.  The boundaries all have flat metrics.  By hypothesis, the injectivity
radii of all points on the boundary are all bounded below so the intrinsic
injectivity radii of boundary tori with respect to the flat metrics will also
be bounded below.  To see that their intrinsic diameters are bounded
above, it suffices to show that their areas are bounded above.  
There are collar neighborhoods of each boundary component with a lower bound on 
their width and the normal curvatures are bounded above.
Thus, if the areas of the boundary were unbounded, the volumes of the collar
neighborhood would be unbounded.   Since the volumes are assumed bounded,
the areas, hence the diameters, are bounded.

The theorems in \cite{Kodani} have the extra hypotheses that the injectivity radius
of all points in the manifold be bounded below, not just boundary points.  Also,
the diameters of the $N_i$ are required to be uniformly bounded above.  However,
the injectivity radius at a point $x$ changes continuously with $x$ and the rate
at which it can go to zero as a function of distance is uniformly bounded 
depending only on the curvature 
(Proposition 8.22 in \cite{GrFrench} or Theorem 8.5 in \cite{GrEnglish}).  
This is often referred to as ``bounded decay of injectivity radius".
It follows that, if the diameters of the $N_i$ are uniformly bounded above, then
the injectivity radius bound on the boundary gives a uniform lower bound to the 
injectivity radius over all of the $N_i$.  The results in \cite{Kodani} apply
directly.

In general, the bounded decay of injectivity radius implies that, if the injectivity
radius at the basepoints of the $N_i$ are bounded below, then, for any fixed distance
$\rho$, the injectivity radius over the neighborhood of radius $\rho$ will be
uniformly bounded below.  The convergence results for manifolds with bounded
diameter give a convergent subsequence for each $\rho$.  The usual diagonal argument
gives a subsequence converging for any fixed $\rho$ which is the definition
of bilipschitz convergence.

Finally, we need to check that the limit manifold is hyperbolic with tubular boundary.
Any interior point in the limit has a neighborhood that is the bilipschitz 
limit of a sequence of
embedded balls in $\H^3$ with fixed radius.  The limit will be isometric to such
a ball so $N_{\infty}$ will be hyperbolic at such a point.  A boundary point
will have a neighborhood that is the bilipschitz limit of a sequence of 
embedded balls on the boundary of  $\H^3_{R_i}$ with fixed radius.  Since
the $R_i$ are bounded below there will be a subsequence which converges to some 
$R$, where possibly $R = \infty$.  The limit neighborhood will be isometric 
to such a ball in $\H^3_R$ so $N_{\infty}$ will have tubular boundary.
\end{proof}

\begin{remark} Although we have based our proof of Theorem \ref{tubelipconv}
on the very general theorems of Gromov and others, there is a much more
direct proof, following the proof of the compactness result of 
J{\o}rgensen-Thurston in \cite {thnotes} (Theorem 5.11.2).  A sketch of
the argument is as follows: For fixed $\epsilon$, let $N_{[\epsilon, \infty)}$
be the set points where the injectivity radius is at least $\epsilon$.
For sufficiently small $\delta$ (depending only on $\epsilon$), there is a 
covering of $N_{[\epsilon, \infty)}$
by embedded balls of radius $\delta$ so that the balls of radius
$\delta/2$ with the same centers are disjoint.  If $N$ is a hyperbolic
$3$-manifold with tubular boundary with $\vol (N) < V$, then the number
of such disjoint balls is bounded in terms of $V$.  Thus, there are
finitely many intersection patterns of the larger balls that cover, 
and the hyperbolic structures on 
$N_{[\epsilon, \infty)}$ are completely determined by the relative positions
of the centers of the balls.  The space of choices of such relative 
positions is compact.
On the other hand, an application of the Margulis lemma, extended to allow
tubular boundary, implies that,
for sufficiently small $\epsilon$ (universal over all hyperbolic $3$-manifolds), 
the regions where the injectivity radius is less than 
$\epsilon$ is a finite disjoint union of tubular neighborhoods of
short geodesics or of cusps.  In the discussion above of canonically 
filling in tubular boundaries, we showed that these regions are determined
isometrically by their boundary data.  This implies 
Theorem \ref{tubelipconv}.

Rather than filling in the details of this argument, we have chosen
to base our proof on published results.  However, some readers may
find this argument clearer.
\end{remark}

Theorem \ref{tubelipconv} allows for the possibility that, even if
all the hyperbolic manifolds $N_i$ are diffeomorphic, the limiting
manifold $N_{\infty}$ may not be.  For this to occur the diameters
must go to infinity.  If this were to occur, then a priori a portion
of the approximating manifolds might be pushed an infinite distance
from the basepoint and be lost in the limit.  This is a familiar
occurrence for hyperbolic surfaces where the length of a geodesic
can go to zero, creating a new cusp and a new diffeomorphism type.

We prove below that this is not possible for sequences of $3$-manifolds 
with tubular boundary having bounded volume and a lower bound for injectivity 
radius at boundary points. 
First we need to establish the fact that the ends of a finite volume
hyperbolic $3$-manifold with tubular boundary have the same structure
as those of a complete, finite volume hyperbolic $3$-manifold.
They are cusp neighborhoods, diffeomorphic to $T^2 \times (0,\infty)$,
formed by dividing out a horoball by a discrete $\Z \oplus \Z$ lattice.
The usual proof that this is the structure of the ends of a
complete, finite volume hyperbolic $3$-manifold uses a refined version
of the Margulis lemma and relies on discreteness of the holonomy group.
The holonomy groups of hyperbolic $3$-manifolds with tubular boundary
are usually not discrete so the proof doesn't immediately apply.
It is possible to give a direct geometric proof for the case with
tubular boundary as in Gromov's extension of the Margulis lemma
(\cite [Proposition 8.51] {GrFrench}).  Instead we use known results about the ends
of finite volume manifolds with pinched negative curvature, due to
Eberlein.

To apply these results we first prove the following lemma:

\begin{lemma} \label{tubemannegcurv} The metric on a hyperbolic
manifold $N$ with tubular boundary can be extended to a complete metric
with pinched negative curvature on a manifold $\M$ diffeomorphic
to the interior of $N$.  $N$ embeds isometrically in $\M$
in this metric and the volume of its complement $\M - N$ is finite.
\end{lemma}

\begin{proof}
The idea of the proof is simply to attach to each component of
the tubular boundary a space diffeomorphic to $T^2 \times (-\infty,0]$,
with $T^2 \times 0$ attached to the boundary.  The result is clearly
diffeomorphic to the interior of $N$.  Furthermore, the metric
on each of the $T^2 \times (-\infty,0]$ pieces will have pinched
negative curvature, finite volume, and agree with the metric on
$N$ in a neighborhood of the tubular boundary.

If $R=\infty$ for a boundary component of $N$, then, as discussed
above, the canonical extension of the boundary metric results in
a finite volume cusp.  In this case, the attached piece has
{\it constant} curvature $-1$.

If $R$ is finite, we use the fact that the metric in a neighborhood
of the tubular boundary is induced from the metric (\ref{cylmetric})
in a neighborhood of $r=R$ by dividing out by the action of
a $(\theta, \zeta)$ lattice.  We alter the metric, keeping it of the form
\begin{eqnarray} dr^2 ~+~ f(r)^2 \,d\theta^2 ~+~ g(r)^2 \, d\zeta^2,
\label{cylnegcurv}
\end{eqnarray}
where $f(r), g(r)$ are defined on $(-\infty, R]$  and agree with $\sinh r,
\cosh r$, respectively near $r=R$.  Furthermore, we want
$f(r), f'(r), f''(r), g(r), g'(r), g''(r)$ to be positive on $(-\infty, R]$.
Such a metric is complete and has negative curvature.  From the explicit 
formulae for the curvatures, it is not hard to see that the sectional
curvatures can be pinched between two negative constants.
(See \cite{agol2} or \cite{kojima} for details of the curvature computation.)

Since the functions $f(r), g(r)$ depend only on $r$, such a metric is
invariant under any $(\theta, \zeta)$ lattice so it descends to a pinched
negatively curved metric on $T^2 \times (-\infty,R]$ which can be attached
to the boundary of $N$.  Further choosing the functions so
that $\int_{-\infty}^R f(r)\, g(r) \, dr < \infty$ ensures that the volume
will be finite.
\end{proof}

\begin{proposition}\label{standardends} Each end of a complete, finite volume 
hyperbolic $3$-manifold with tubular boundary is diffeomorphic to 
$T^2 \times (0, \infty)$ and is isometric to a horoball in $\H^3$ divided 
out by a parabolic $\Z \oplus \Z$ lattice.
\end{proposition}

\begin{proof}  In \cite{Eb} it is proved that for complete, finite volume
$n$-manifolds with pinched negative curvature, there will be a finite number
of ends, each of the form
$W \times (0, \infty)$ where $W$ is an $(n-1)$-manifold with virtually
nilpotent fundamental group.  Since our manifolds are orientable and
$3$-dimensional, $W$ is an orientable surface and the only possibility
is a torus. It is further shown that the end is isometric to a horoball
divided out by a parabolic lattice isomorphic to the fundamental group 
of $W$.  In the general negatively curved context, horoballs are defined 
in terms of Busemann functions.  However, since the ends of the negatively 
curved manifold constructed in Lemma \ref{tubemannegcurv} that come from 
the original hyperbolic manifold with tubular boundary all have constant 
curvature, a horoball sufficiently far out in the end defined by a Busemann 
function will agree with the usual definition in hyperbolic geometry.
\end{proof}

We are now in a position to prove a compactness result for the
set of hyperbolic structures with tubular boundary on a fixed
compact $3$-manifold.

\begin{theorem} \label{tubelipcpt} The set of hyperbolic structures
with tubular boundary on a fixed compact $3$-manifold $N$ with volumes bounded
above and injectivity radius on the boundary bounded below is compact
in the bilipschitz topology.  In other words, suppose that $N_i$ is
a sequence of hyperbolic manifolds with tubular boundary, all diffeomorphic
to $N$.  Assume there are constants $c, V >0$ such that, for all $i$,
$inj(x, N_i) \ge c$ for all $x \in \bd N_i$ and $\vol(N_i)\leq V$.
Then there is a subsequence which converges in the bilipschitz topology
to a hyperbolic structure on $N$ with tubular boundary.
\end{theorem}

\begin{proof}
By Theorem \ref{tubelipconv}, it suffices to show that the diameters
of the $N_i$ are uniformly bounded.
Choose a basepoint $x_i \in \bd N_i$ for all $i$.
Again, by Theorem \ref{tubelipconv}, there will always be a subsequence 
of $(N_i,x_i)$ with a limit $(N_{\infty},x_{\infty})$ in the pointed 
bilipschitz topology which is again a pointed hyperbolic $3$-manifold 
with tubular boundary.

Suppose that the diameters of the $N_i$ are not bounded above.
By definition of convergence in the bilipschitz topology the
limit will be non-compact and will have finite volume.
It will have at least one end, and each end is a cusp with
a horospherical Euclidean torus as cross section by
Proposition \ref{standardends}.

Convergence in the  bilipschitz topology further implies that
we get a sequence of bilipschitz maps
of larger and larger radius neighborhoods of $x_{\infty} \in N_{\infty}$
into $N_i$ which, for any fixed radius, are becoming
arbitrarily close to an isometry onto their images.  For a sufficiently
large radius, the topology of these neighborhoods will be constant and
equal to a manifold $W$ with torus boundary components whose interior
is diffeomorphic to $N_{\infty}$.  We identify the $N_i$ with $N$ and
the large radius neighborhoods with $W$ and consider the bilipschitz maps
as maps $\phi_i: W \to N$.  Under the identification of the interior of $W$
with $N_\infty$ certain of the boundary tori of $W$ correspond to cusps
of $N_\infty$.  We will refer to these tori as ``cusp tori".

The hyperbolic structures, $N_i$, on $N$ induce holonomy representations
$\rho_i: \pi_1 N \to G$, where $G$ is the group of isometries of $\H^3$.
The representations are well-defined up to conjugation by elements in $G$.
Similarly, the hyperbolic structure $N_\infty$, viewed as a structure on
the interior of $W$, induces a representation $\rho: \pi_1 W \to G$.  
The fact that the bilipschitz maps converge on compact sets implies the 
convergence of the holonomy representations of any finite set of group elements, 
at least after conjugating the representations.
Since $\pi_1 W$ is finitely generated, this implies that, perhaps after
conjugating the $\rho_i$ by elements of $G$, we obtain
\begin{eqnarray}
\rho_i \circ (\phi_i)_{*} \to \rho.
\label{holconv}
\end{eqnarray}

By Proposition \ref{standardends} the fundamental group of the 
torus cross-sections of the cusp ends of $N_\infty$ inject into the 
fundamental group of $N_\infty$.  Since $N_\infty$ is diffeomorphic to the 
interior of $W$, it follows that the fundamental group of each cusp torus 
of $W$ injects into the fundamental group of $W$.  Choose any cusp torus 
and denote it by $T$.  We wish to show that, for $i$ sufficiently large,
the fundamental group of $T$ must inject under $(\phi_i)_{*}$ into
the fundamental group of $N$.  Furthermore, $T$ will not be
peripheral in $N$.  This will contradict the fact that $N$ is atoroidal,
implying that the diameters of the $N_i$ must have been uniformly bounded above.

For each value of $i$ we denote by $W_i$ the homeomorphic
image of $W$ in $N$ under $\phi_i$ and by $T_i$ the homeomorphic
image of $T$.
Suppose, for some $i$,  the torus $T_i \subset W_i$ is
compressible in $N$.  Since $N$ is irreducible, the torus must either
bound a solid torus outside $W_i$ or be contained in a 3-ball in $N$.
For any element $\gamma \in \pi_1 T$ we have
$\rho_i \circ (\phi_i)_{*}(\gamma) \to \rho(\gamma)$.
Since, for any non-trivial $\gamma$, $\rho(\gamma)$ is a non-trivial
parabolic element, this implies that $\rho_i \circ (\phi_i)_{*}(\gamma)$ 
is non-trivial for sufficiently large $i$.  Hence, $\pi_1 T$ at least 
maps non-trivially under $(\phi_i)_{*}$. Therefore, no cusp torus is 
contained in a 3-ball so all the cusp tori must bound solid tori outside
$W_i$.  Since this is true for all of the cusp tori in $W$, it follows
that, for all sufficiently large $i$, adding $N - W_i$ to $W_i \subset N$
corresponds to obtaining $N$ by Dehn filling on $W$.

Let $\gamma_i$ denote a curve on a cusp torus $T$ of $W$ which
bounds a disk when mapped into $N$ by $\phi_i$.  As above, for
any fixed non-trivial element $\gamma \in \pi_1 T$, 
$\rho_i \circ (\phi_i)_{*}(\gamma)$ will be non-trivial for
sufficiently large $i$ (where ``sufficiently large" generally depends
on $\gamma$).  Since $(\phi_i)_{*} (\gamma_i) = e$,
its holonomy representation is trivial.   Thus, $\gamma_i$ can represent 
a fixed element of $\pi_1 T$ for only {\em finitely many} values of $i$.  
Since this argument holds for each cusp torus, it implies that $N$ 
can be obtained by Dehn fillings on $W$ using
infinitely many distinct filling curves on each cusp torus.  We will show that 
this is impossible by Thurston's theory of hyperbolic Dehn surgery.

First, note that, since $N$ has a complete metric of pinched negative
curvature, it is irreducible and atoroidal (\cite{Eb}).  It is the interior of the 
compact manifold $W$ with non-empty boundary which is therefore Haken.
By Thurston's Geometrization Theorem for Haken manifolds (\cite{th-unif},
\cite{otal}, \cite{otal2}, \cite{kap}), $N$ supports a complete, finite volume metric 
of constant negative curvature.
Thurston's hyperbolic Dehn surgery theorem says that, when considering 
all possible Dehn fillings of such a $3$-manifold, for all but
finitely many choices of filling curve on each cusp torus, 
the result is hyperbolic.  Thus, for $i$ sufficiently large, all the manifolds 
obtained above by Dehn filling $W$ are hyperbolic.
Furthermore, they have volumes converging from below to the volume
of the complete hyperbolic structure on $N$.  
But, since the resulting $3$-manifold is always diffeomorphic to $N$ and 
the hyperbolic volume of $N$ is a topological invariant, this is
a contradiction.
\end{proof}

\begin{remark}The above result generalizes to the case when $N$
has cusps.  To do this, one shows, (using, for example, the packing
results of the next section), that it is possible to remove
neighborhoods of the cusps in such a way that the injectivity radii of
the new boundary components of the resulting compact hyperbolic manifold
with tubular boundary are also bounded below.
\end{remark}

We are now in a position to prove our main convergence result,
referred to in the introduction as Theorem \ref{thm2}.

\begin{theorem} \label{newthm2}  Let $\N_t$, $t\in [0,t_\infty)$, be 
a smooth path of closed hyperbolic cone-manifold structures
on $(\N, \Sigma)$ with cone angle $\alpha_t$ along the singular locus $\Sigma$.
Suppose $\alpha_t \to \alpha \geq 0$ as $t \to t_\infty$, that the volumes of 
the $\N_t$ are bounded above by $V_0$, and that there is a positive constant
$R_0$ such that there is an embedded tube of radius at least $R_0$ around
$\Sigma$ for all $t$.  Then the path extends continuously to $t=t_\infty$ 
so that as $t \to t_\infty$, $\N_t$ converges in the bilipschitz topology 
to a cone-manifold structure $\N_\infty$ on $\N$ with cone angles 
$\alpha$ along $\Sigma$.
\end{theorem}

\begin{proof} Removing disjoint tubular neighborhoods of the singular
locus, we obtain a smooth path of hyperbolic manifolds $N_t$ with tubular
boundary, with all the $N_t$ diffeomorphic to a fixed compact $3$-manifold $N$.  
The volumes are bounded above since they are smaller
than the volumes of the cone-manifolds $\N_t$ which are bounded
above by hypothesis.

To apply Theorem \ref{tubelipcpt} we need to show that
there is a lower bound to the injectivity radii on the boundary
of the $N_t$. That will imply that there is a subsequence of the $N_t$
converging to a hyperbolic manifold $N_{\infty}$ with tubular boundary.
The boundary can then be filled in canonically to obtain a hyperbolic 
cone-manifold $\N_{\infty}$.

By definition, the injectivity radius at a boundary point is less
than its distance to the singular locus in the corresponding
hyperbolic cone-manifold structure.  Similarly, it is less than
its distance to any other boundary components besides the one it is on.
Since the tube radii of the hyperbolic cone-manifolds are bounded below by 
hypothesis, the tubular neighborhoods that are removed can be chosen so that 
both the distance to the singular locus and to other boundary components 
are bounded below.  Furthermore, we will see in the next section that, 
on a boundary torus of radius $R$, there is always an embedded ellipse 
with minor axis lengths given by (\ref{minors}).  This implies
that the injectivity radii are bounded below as desired.

Take any sequence $N_{t_j}$, where 
$t_j \in [0,t_\infty), t_j \to t_\infty$.
We can apply Theorem \ref{tubelipcpt} to conclude that there
is a subsequence $N_i$ which converges
in the bilipschitz topology to a hyperbolic manifold, $N_\infty$,
with tubular boundary.  It is also diffeomorphic to $N$.  As in
the proof of the previous theorem, the hyperbolic structures
$N_i$ and $N_\infty$ gives rise to holonomy representations
$\rho_i$ and $\rho$ respectively from $\pi_1 N$ to the group $G$ of
isometries of $\H^3$.  Since the diameters of the $N_i$ are
uniformly bounded, convergence in the bilipschitz topology
provides basepoint-preserving bilipschitz homeomorphisms from 
$N_\infty$ to the $N_i$ which, under the identifications of both the 
domain and range with $N$, give basepoint-preserving homeomorphisms 
$\phi_i: N \to N$.  As in the proof of the previous theorem, it is 
possible to choose conjugacy classes of the holonomy representations
so that $\rho_i \circ (\phi_i)_{*} \to \rho$.

Since the $\phi_i: N \to N$ are basepoint-preserving homeomorphisms,
the induced maps $(\phi_i)_{*}$ on $\pi_1 N$ are automorphisms.
We saw in the proof of Theorem \ref{tubelipcpt} that $N$ has a
complete, finite volume hyperbolic metric on its interior.  This 
implies that the outer automorphism group
of  $\pi_1 N$ is finite because Mostow rigidity says that any
outer automorphism is homotopic to an isometry of the complete
finite volume metric on the interior of $N$.  The group of such
isometries is finite.  (See \cite{thnotes} for a more detailed
version of this argument.) Since there are only finitely many
choices for $(\phi_i)_{*}$ up to conjugacy, there is a further
subsequence so that $(\phi_i)_{*}$ is constant and, hence, that
$\rho_i \circ (\phi)_{*} \to \rho$ for a fixed automorphism
$(\phi)_{*}$ of  $\pi_1 N$.
This implies that the holonomy representations $\rho_i$
converge in the representation variety (representations of $\pi_1 N$ to 
$G$ modulo conjugation) to $\rho \circ (\phi)_{*}^{-1}= \hat \rho$.

The hyperbolic structure $N_\infty$ with tubular boundary
has $\hat \rho$ as a holonomy representation.  The boundary
data of the $N_i$ determine the canonical completion to the
hyperbolic cone-manifold structures $\N_i$.  Since these boundary 
data converge to that of $N_\infty$, its canonical completion is
a hyperbolic cone-manifold structure $\N_\infty$ with cone angle 
$\alpha$ along its singular locus $\Sigma$.  
Under the isomorphism $\pi_1 N \cong \pi_1 (\N_\infty - \Sigma)$ 
the holonomy representation of $\N_\infty - \Sigma$ can be identified
with $\hat \rho$.

The local rigidity theorem of \cite{HK1} implies that the hyperbolic
cone-manifold structures on $(\N_\infty, \Sigma)$ with a fixed cone angle 
(with angle at most $2\pi$) are isolated.  The above analysis applies to any 
convergent subsequence of the $\N_t$.  If we view the path $\N_t$ as
a path $\rho_t$ in the representation variety, this implies that any
accumulation point of $\rho_t$ as $t \to t_\infty$ corresponds 
to a hyperbolic cone-manifold structure on $(\N_\infty, \Sigma)$ 
with cone angle $\alpha$.  Since these are isolated and the set of 
accumulation points is connected, there can be only a single 
accumulation point.  It follows that the path $\rho_t$ extends 
continuously to $t_\infty$ and that the $\N_t$ converge in the 
bilipschitz topology to $\N_\infty$.
\end{proof}

\section{A packing argument}\label{pack}

Let $M$ be a $3$-dimensional hyperbolic cone-manifold with 
a link $\Sigma$ as singular locus.
Let $R$ be the radius of the maximal embedded tube in $M$
around $\Sigma$ and denote this tube by $U_R$. 
If $\Sigma$ has multiple cusps, this is to be interpreted as meaning that the
radii of the tubes around all of the components are the same, equal to $R$.
In this section we will find lower bounds for the area of each component 
of the boundary of $U_R$
via a packing argument analogous to the usual horoball packing
arguments for non-singular cusped hyperbolic $3$-manifolds
(cf. \cite{Me}, \cite{adams}). For non-singular hyperbolic
3-manifolds, similar tube packing arguments are used in \cite{gmm}.

Denote by $\tilde X$
the universal cover of $X= M - \Sigma$, equipped
with  the lift of the metric on $X$.  The developing map $\tilde X \to
\H^3$ can be extended by completion to the lifts of the 
singular locus, giving a map $D: \hat M \to \H^3$ where 
$\hat M$ is the metric completion of $\tilde X$. Further the covering
projection $\tilde X \to X$ extends by completion to a map $p: \hat M \to M$.
($\hat M$ can be regarded as the universal branched covering of $M$,
 branched over $\Sigma$.)

Choose a component, $\Sigma_0$, of the lift of a component of the 
singular locus to $\hat M$.  Under the developing map $\Sigma_0$
maps to a geodesic, $g$, in $\H^3$. The universal cover of 
$\H^3 - g$ can be completed by adding a geodesic, $\hat g$, which projects 
to $g$ in $\H^3$. (This can be thought of
as the infinite cyclic branched cover of $\H^3$ branched
over the geodesic $g$.) Let $\hat \H^3$ denote this completion and let 
$\hat U_r$ denote the neighborhood of radius $r$ about $\hat g$ in $\hat \H^3$.
Then for each $r < R$, the $r$-neighborhood $U_r$ of each component of 
$\Sigma$ in $M$ is isometric to the quotient of
$\hat U_r$ by a discrete group $\Gamma \cong \Z \oplus \Z$ of isometries
of $\hat \H^3$ preserving the axis $\hat g$.

We can also regard $\hat \H^3$ as the
``normal bundle'' to $\Sigma_0$ in $\hat M$ and there
is an exponential map $E: \hat U_{2R} \to \hat M$ defined
by extending geodesics orthogonally from $\Sigma_0$.
This gives a isometric embedding from $\hat U_{2R}$
onto the neighborhood of radius $2R$ about $\Sigma_0$ in
$\hat M$. 

Because $R$ is the maximal tube radius, there 
is a geodesic arc $\tau$ of length $2 R$ in $M$ going from $\Sigma$ to 
itself which is perpendicular to $\Sigma$ at both endpoints.
It is a shortest geodesic arc from $\Sigma$ to itself not entirely 
contained in $\Sigma$.  The radius $R$ tube around $\Sigma$, $U_R$, 
has a self-tangency at the midpoint of $\tau$.  Now consider all the 
lifts to $\hat M$ of arcs of length $2 R$ from $\Sigma$ to itself, 
beginning at 
$\Sigma_0$. They end at points $q_i$ lying on other lifts of $\Sigma$.  
These points can be identified, via the inverse of the exponential 
map $E$ with points, also denoted by $q_i$, in $\hat U_{2R}$.

\begin{lemma} Let $\{ q_i \}$ be the set of all endpoints 
of such lifts of arcs of length $2 R$ from $\Sigma$ to itself. 
Then the distance between $q_i$ and $q_j$ in $\hat \H^3$ satisfies 
$d(q_i,q_j) \ge 2R$ for all $i \ne j$.
\end{lemma}

\begin{proof} Consider two points $q_i,q_j$ with $i\ne j$. 
These lie on the boundary of the set $\hat U_{2R}$
in $\hat \H^3$, which is convex since the distance to a geodesic is
a convex function.  Thus, 
the shortest geodesic $\gamma$ from $q_i$ to $q_j$ in $\hat \H^3$ lies
inside $\hat U_{2R}$. Composing $\gamma$ with the exponential
map  $E: \hat U_{2R} \to \hat M$ and the (branched) covering 
projection $p: \hat M \to M$ gives a geodesic  $\bar \gamma$ in 
$M$ joining $\Sigma$ to itself.  Since $\bar\gamma$ is not entirely 
contained in $\Sigma$ it has length at least $2R$. 
Hence $d(q_i,q_j) \ge 2R$.
\end{proof}

%\medskip

For each $i$, let $B_i$ denote the ball in $\hat \H^3$ of radius 
$R$ about $q_i$.  We project the balls for all the 
$q_i$ orthogonally onto the surface 
$\bd \hat U_R$ in $\hat \H^3$ at radius $R$ from the singular set. 
The fact that the balls $B_i$ are disjoint implies that their 
projections $P_i$ are also disjoint.  This follows easily
from the facts that the centers of the 
$B_i$ all have the same radial co-ordinate and all of the 
$B_i$ have the same radius. 

Next we will estimate the area of each $P_i$ and use this to 
estimate the area of $T_R$.
But first we prove some preliminary geometric results. 

Let $(r,\theta,\z)$ denote hyperbolic
cylindrical coordinates on $\H^3$ about a geodesic $g$. These
can also be regarded as cylindrical co-ordinates on $\hat \H^3$ 
about the geodesic $\hat g$ covering $g$, but the angle
$\theta$ is no longer measured modulo $2 \pi$, but
rather as a real number. 

\begin{lemma}\label{cyldist} 
The distance $d$ between two points $p_1,p_2$ in $\hat\H^3$ with cylindrical coordinates
$(r_1,\theta_1,\z_1)$ and $(r_2,\theta_2,\z_2)$ with $|\theta_1 - \theta_2 |\le \pi$
is given by
$$\cosh d = \cosh(\z_1-\z_2) \cosh{r_1} \cosh{r_2} -
\cos(\theta_1-\theta_2) \sinh{r_1} \sinh{r_2} .$$
\end{lemma}

\begin{proof} See \cite{gmm}, Lemma 2.1. \end{proof}

%\smallskip

We now study the projection of a ball onto a hyperbolic
cylinder.

\begin{lemma} \label{projectball} Consider a ball of radius $d$ centered at the
point with cylindrical coordinates $(r,\theta,\z) =
(r_0,0,0)$ with $d < r_0$.
The projection of this ball to the $(\theta,\z)$-plane
has equation
$$\sinh^2 \z \cosh^2 r_0 + \sin^2\theta \sinh^2 r_0 \le \sinh^2d.$$
\end{lemma}

\begin{proof} From the distance formula in cylindrical coordinates (Lemma \ref {cyldist}),
the ball has equation
$$\cosh \z \cosh {r_0} \cosh r -\cos\theta \sinh {r_0} \sinh
r \le \cosh d.$$
Writing $\cosh r$ and $\sinh r$ as exponentials gives
$$\cosh \z \cosh {r_0} (e^r + e^{-r}) - 
\cos\theta \sinh {r_0} (e^r - e^{-r}) - 2 \cosh d \le 0$$
or
$$e^{2r}(\cosh \z \cosh r_0 - \cos\theta \sinh r_0)
-2e^r \cosh d + (\cosh \z \cosh r_0 + \cos\theta \sinh r_0) \le 0.$$
Given $(\theta,\z)$ this quadratic for $e^r$ has a real solution
if and only if the discriminant is non-negative, i.e.
$$(2\cosh d)^2 - 4(\cosh \z \cosh r_0 - \cos\theta \sinh r_0)
(\cosh \z \cosh r_0 + \cos\theta \sinh r_0) \ge 0,$$
or
$$\cosh^2 \z \cosh^2 r_0 - \cos^2 \theta \sinh^2 r_0 \le \cosh^2
d.$$
Rewriting this, using $\cosh^2 \z = \sinh^2 \z + 1$ and $\cos^2\theta
= 1 - \sin^2 \theta$, we
have
$$\sinh^2 \z \cosh^2 r_0 + \sin^2\theta \sinh^2 r_0 \le \sinh^2d.$$
\end{proof}

Each ball $B_i$ has radius $R$ and its center is at distance $2 R$ from the
geodesic $\hat g$ in $\hat \H^3$.   We choose co-ordinates so that the 
center of a ball $B_i$ has co-ordinates $(r,\theta,\z) = (2 R, 0,0)$.  
From Lemma \ref{projectball}, the projection
of $B_i$ onto the $(\theta, \z)$-plane satisfies the equation:
\begin{eqnarray}
f(\z,\theta) = \sinh^2 \z \cosh^2 2R + \sin^2\theta \sinh^2 2R \le \sinh^2R.
\label{projeq}\end{eqnarray}

Ignoring the self-tangencies, the boundary of the maximal tube, 
$U_R$, in the cone-manifold is a torus $T_R$ with an induced 
Euclidean structure.  The Euclidean structure is induced
from the set of points in $\hat \H^3$ at distance $R$ from
$\hat g$ modulo the group $\Gamma$.  Since the projection of each
$B_i$ onto the $(\theta, \z)$-plane is disjoint from its translates
under $\Gamma$, the corresponding set $P_i$ with radial co-ordinate
$R$ is disjoint from its translates.  This implies that
it embeds in $T_R$ under the quotient map from the action
of $\Gamma$. Further the collection of $P_i$ contains
at least two distinct $\Gamma$-orbits if $\Sigma$ consists
of a single component.

Next we estimate the area of each $P_i$ and use this
to estimate the area of $T_R$.

\begin{theorem} \label{areathm}
The area of the torus $T_R$ at distance $R$ from $\Sigma$
satisfies
\begin{eqnarray}
{\rm area}(T_R) \ge  3.3957 {\sinh^2 R \over \cosh(2R)}.
\label{areabound}\end{eqnarray}
if $\Sigma$ is connected.  If $\Sigma$ has multiple components,
then, for each component, the lower bound for the area of the torus $T_R$ is
half as large: 
\begin{eqnarray}
{\rm area}(T_R) \ge  1.6978 {\sinh^2 R \over \cosh(2R)}.
\label{multipleareabound}\end{eqnarray}

\end{theorem}

\begin{proof} 
Equation (\ref{projeq})  gives us bounds on $\z$ and $\theta$: 
$$|\sinh \z| \le  {\sinh(R) \over \cosh(2R)} {\rm~~and~~}
|\sin\theta| \le {\sinh(R) \over \sinh(2R)}.$$ 

Now $${\sinh R \over \cosh(2R)} = {s \over 1+ 2s^2}$$ where $s = \sinh(R)$.
By the arithmetic-geometric mean inequality
$\sqrt{2} s = \sqrt{2s^2} \le {1+ 2s^2 \over 2}$, hence ${s \over 1+ 2s^2}
\le {1\over 2\sqrt{2}}$ for all $s \ge 0$, with equality attained exactly
when $1=2s^2$, i.e. $R=R_0$ where $\sinh(R_0) = {1\over \sqrt{2}}$. 

So for such $\z$
we have  $|\sinh \z| \le {1\over 2\sqrt{2}}$. 
Since $\sinh \z$ is a convex function for positive values of $\z$, we obtain
$\left|{\sinh \z \over \z}\right| \le \S$ where
$$\S = {{1\over 2\sqrt{2}} \over {\arcsinh({1\over2\sqrt{2}})}}
\approx {1 \over 0.980258}.$$
Since $|{\sin\theta \over \theta}|\le 1$, we deduce that
$$f(\z,\theta) \le (\S \z)^2 \cosh^2(2R) + \theta^2 \sinh^2(2R).$$ Thus the projected ball
defined by equation (\ref{projeq})
contains the region
$$(\S \z)^2 \cosh^2(2R) + \theta^2 \sinh^2(2R) \le \sinh^2 R$$
or
\begin{eqnarray}\left({\S\cosh(2R) \over \cosh R \sinh R}\right)^2(\z \cosh R)^2 + 
\left({\sinh(2R) \over \sinh^2 R}\right)^2
(\theta\sinh R)^2 \le 1. \label{elleq}\end{eqnarray}  
Since $\z \cosh R$ and $\theta\sinh R$ are Euclidean coordinates on the torus at radius $R$,
equation $(\ref{elleq})$ describes an ellipse with semi-major axes 
\begin{eqnarray}\label{minors}
a = {\cosh R \sinh R\over \S \cosh(2R)} {\rm\quad and \quad} b = {\sinh^2 R\over
\sinh(2R)}\end{eqnarray} 
and area $$\pi a b = {\pi \sinh^2 R\over 2 \S \cosh(2R)}.$$ The axes of all of 
the ellipses are parallel to the $\theta$ and $\z$ axes. By an area preserving affine
transformation of the torus, we can arrange that all the inscribed ellipses simultaneously
become circles of the same radius. It follows that the packing density of the
ellipses is at most the maximum packing density of circles, 
namely $\pi \over {2 \sqrt{3}}$.

Furthermore, if $\Sigma$ is {\em connected},  the torus $T_R$ at radius $R$ contains 
at least {\em two} disjoint ellipses, so its area satisfies:
$${\rm area}(T_R) \ge {2 \sqrt{3} \over \pi} 2 \pi ab =\sqrt{3} a b
= {2 \sqrt{3} \sinh^2 R \over  \S \cosh(2R)},$$
so $${\rm area}(T_R) 
\ge 3.3957 {\sinh^2 R \over \cosh(2R)}.$$
If $\Sigma$ has multiple components, then, for each component, 
the lower bound for ${\rm area}(T_R)$ is half as large.
\end{proof}

\section{Controlling the tube radius}\label{tuberad}

In this section we will use the information derived in 
sections \ref{harmdef} and \ref{pack} to control the change 
in the radius of the maximal
embedded tube around the singular locus.  This will allow us to
complete the proof of Theorem \ref{thm1}. 
Finally we combine this with Theorem \ref{thm2} to prove Theorems \ref{thm3}
and \ref{thm4}.

Rather than studying the tube radius directly, we will derive
information about it by studying the geometry of the torus
on the boundary of the maximal tube.  The boundary torus
has an intrinsic flat metric.  We will denote by $\m$ the
length in this metric of the geodesic in the homotopy class
of the meridian.  The height of the maximal annulus with the
meridian as its core will be denoted by $\h$.  Thus, the area of
the torus, denoted by $A$, will equal $\m \h$.  If the radius
of the tube is $R$, then $\m$, $\h$ and $A$ are related to the
cone angle $\alpha$ and the length $\l$ of the singular
locus by the formulae:
\begin{eqnarray*} \m &=& \alpha \sinh R,\\
 \h &=& \l \cosh R,\\
 A &=& \alpha\l \sinh R \cosh R.\end{eqnarray*}

Theorem \ref{areathm} implies that the area $A$ of the flat torus 
satisfies 
$$A \geq 3.3957 {\sinh^2 R \over \cosh\, (2R)}.$$  Dividing
by $\sinh R \cosh R$ provides the following key estimate. 

\begin{corollary} 
Suppose the singular set $\Sigma$ has length $\l$ and cone angle $\alpha$. Then
the radius $R$ of 
a maximal embedded tube about $\Sigma$  satisfies
\begin{eqnarray} \alpha \l \ge h(R)  
= 3.3957 {\tanh R \over \cosh \,(2R)}. \label{humpeq}\end{eqnarray}
\end{corollary}

%\medskip
\begin{remark*} In the case of a closed geodesic in a non-singular 
hyperbolic $3$-manifold we have $\alpha = 2\pi$, and this gives
$$\l \ge 0.5404 {\tanh R \over \cosh \,(2R)}.$$
This seems to be very close to the estimate given in Proposition
3.1 of \cite{gmm}. 
\end{remark*}

\smallskip

The qualitative behavior of the function $h(r) = 3.3957 {\tanh r \over \cosh \,
(2r)}$, 
whose graph is pictured below,
is very important and will influence the form of all of our arguments.

\bigskip
\centerline{\epsfbox{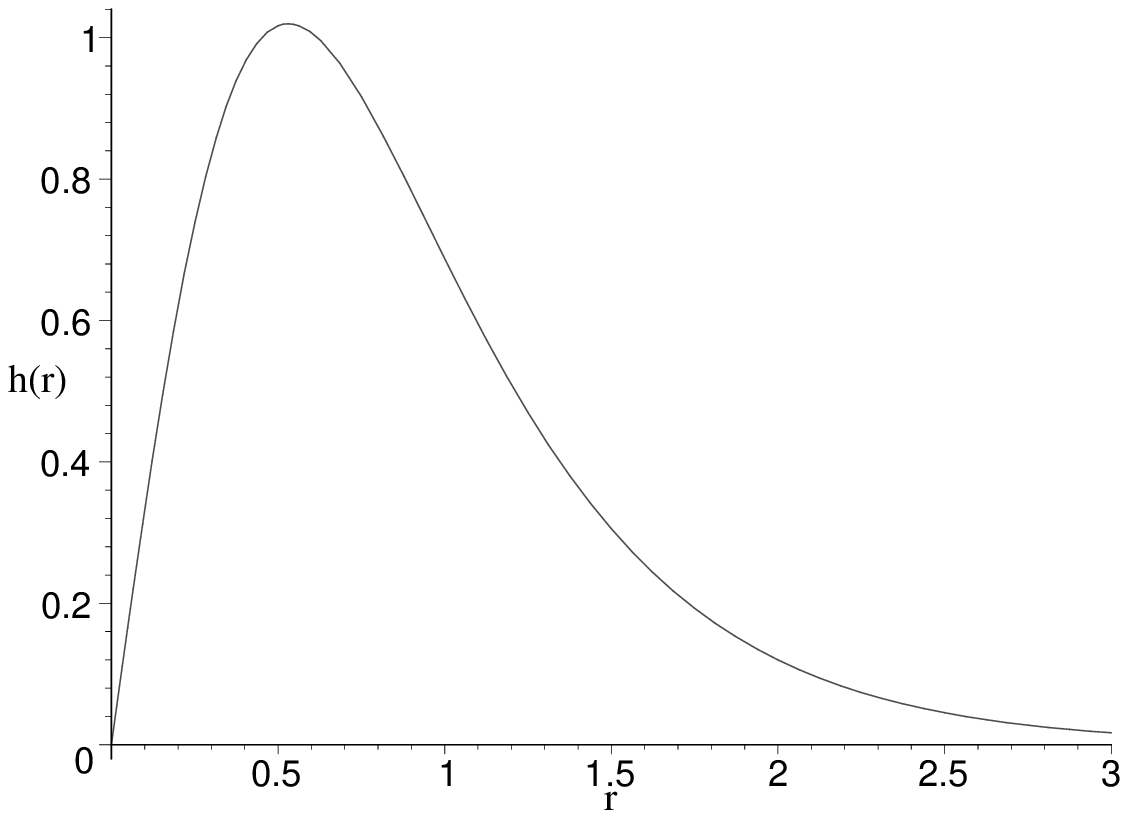}}
%
%\bigskip

The inequality (\ref{humpeq}) implies that, for a given tube radius,
there is a lower bound
to the product, $\alpha \l $.  Hence, for a given tube radius and cone angle,
a lower bound to the length of the core curve.  Instead, we would like to bound
the tube radius in terms of $\alpha \l $.  This is, in fact, not literally
possible and is reflected by the graph of $h(r)$ as it drops down
to $0$ as $r \to 0$.  

However, note that $h(r)$ appears to have a single maximum near $r=.5$
(more precise values are given below) and to be strictly decreasing for
values of $r$ larger than this.  In particular, it appears to be invertible
for such values of $r$.  Thus, if the tube radius $R$ is known to 
be larger than this value and if $\alpha \l $ is smaller than the 
maximum value of $h(r)$, then the value $h^{-1}(\alpha \l)$ of the 
inverse function will provide a further lower bound for $R$. This lower
bound goes to infinity as $\alpha \l $ goes to zero.

In our situation, we will be starting with a complete structure, for which
the tube radius is infinite and $\alpha = \l = 0$.  In particular, 
as we try to increase the cone angle, we begin with values of 
the tube radius and $\alpha \l $ 
for which the inverse of the function $h(r)$ provides a lower bound to the
tube radius.  As long as the value of $\alpha \l $ remains below the
maximum value of $h(r)$, the tube radius is bounded below and the results
of section 3 imply that there can be no degeneration.

\smallskip

The goal of this section is to provide initial conditions on the surgery
curve that will guarantee that $\alpha \l $ remains below this maximum
value until the cone angle reaches $2 \pi$.

%\medskip

\begin{remark*} For smooth structures, i.e. when $\alpha = 2 \pi$,  the
results of \cite{meyerhoff}
imply that, for sufficiently short geodesics, there is a lower bound to
the tube radius.  This result uses J{\o}rgensen's inequality, which has no 
literal analogue for cone-manifolds.
To see that there is no such lower bound for the tube radius around short
core curves in a general cone-manifold, one can consider the figure eight 
knot complement and choose the standard meridian as the surgery curve. 
As the cone angle increases, the length of the core geodesic increases for 
a while (enough for $\alpha\l$ to become larger than the maximum of $h(r)$), but 
then goes to $0$ as the cone angle approaches 
$\alpha = {{2 \pi}\over{3}}$.  In fact, the hyperbolic structures 
degenerate in such a way, that, if they are rescaled to have 
volume $1$, they converge to a Euclidean orbifold
at  $\alpha = {{2 \pi}\over{3}}$ .
\end{remark*}

The following lemma shows that the qualitative behavior of the function 
$h(r)$ which was presumed in the previous discussion is as desired.  It
also provides an accurate value for the maximum of $h(r)$ and for the value
of $r$ at which it is attained.

%\medskip

\begin{lemma}\label{hinverse}
The function $h(r)$
is a decreasing function of $r$ for $r \ge 0.531$ with an 
inverse $h^{-1}(a)$ defined for $0\le a \le h_{max} = h(0.531) \approx 1.019675$
such that $h^{-1}(a)=r$ if and only if $h(r)=a$ and $r\ge 0.531$.
\end{lemma}

\begin{proof}
Writing the function $h$ in terms of $\z =\tanh r$, we have
$$h(r) = 3.3957 ~{\tanh r} ~ { \cosh^2 r - \sinh^2 r \over \cosh^2 r + \sinh^2 r}
= 3.3957 ~ {\z (1-\z^2) \over 1+\z^2}.$$
If we put $f(\z)=\displaystyle{\z (1-\z^2) \over 1+\z^2}$, then
$f'(\z) = \displaystyle{1-4\z^2 - \z^4 \over (1+\z^2)^2}$. Hence
$f(\z)$ has a unique maximum for $0<\z<1$ when
$1-4\z^2 - \z^4 =0$, or $\z^2 = \sqrt{5}-2$. Then
$\z\approx 0.485868$, $r =\arctanh(\z) \approx 0.5306375$ and
$h(r)=3.3957 \, f(\z) \approx 1.0196755$. The result follows immediately.
\end{proof}

From this lemma we deduce that the estimate (\ref{humpeq}) 
gives us a lower bound 
for the tube radius in terms of $\alpha \l$: 
%\bigskip

\begin{proposition}\label{tuberadbound} The tube radius $R$ satisfies
\begin{eqnarray} R \ge h^{-1} (\alpha \l) {\rm~when~}\alpha \l \le h_{max} 
\approx 1.019675 {\rm~and~}
R \ge 0.531.
\label{Rbound}\end{eqnarray} 
\end{proposition} 

Together with the non-degeneration results of section 
\ref{geomlim} we immediately
have the following theorem:

%\medskip

\begin{theorem}\label{alnondeg} Let $\N_s$ be a smooth family of 
finite volume $3$-dimensional hyperbolic cone-manifolds, with 
cone angles $\alpha_s, 0\le s < 1$, where $\lim_{s\to 1}\alpha_s = \alpha_1$. 
Suppose the tube radius $R$ satisfies $R \ge 0.531$ for $s=0$
and $\alpha_s \l_s \le h_{max}$ holds for all $s$, where $\l_s$ denotes the 
length of the singular geodesic.  If the volumes of the $\N_s$ remain bounded,
then the $\N_s$ converge geometrically to a cone-manifold $\N_1$ with cone 
angle $\alpha_1$.  In particular, this conclusion holds if $\N_0$ is complete 
($\alpha_0 = 0$), 
$\alpha_s$ is increasing and $0 <\alpha_s \l_s \le h_{max}$ for all $s$.
\end{theorem}

\begin{proof} Proposition \ref{tuberadbound}  implies that, 
if the initial tube radius is at least $0.531$, then, since
$h^{-1}(\alpha \l) \geq 0.531$ by definition, the tube radius will remain 
at least $0.531$ as long as $h^{-1}$ is defined.  This will be the case as long 
as $\alpha_s \l_s \le h_{max}$.  The first statement now follows 
immediately from Theorem \ref{thm2}.  
In the special case when $\N_0$ is complete, the tube radius 
is infinite, hence bigger than $0.531$, for $s=0$. From the Schl\"afli 
differential formula (\ref{schlafli}), 
the volume decreases as the cone angle increases. 
Hence the volumes are uniformly bounded throughout the deformation and this
special case follows from the general case.
\end{proof}

\medskip

In light of the above theorem, we would like to find a method to bound the 
quantity $\alpha \l$ from above throughout a deformation.  Since $t = \alpha^2$ 
is our parameter, this amounts to controlling the growth of the core 
length $\l$.
Our estimates from Section \ref{harmdef} provide control of the change 
in $\l$ in 
terms of $\alpha$ provided that the tube radius is bounded below.
Specifically, recall that equation (\ref{dldalpha}) gives
$${d \l \over d \alpha} = {\l \over \alpha} ( 1 + 4\alpha^2 x),$$
and that we have the estimate (\ref{keybound})
$$ {-1 \over \sinh^2(R)} 
	\biggl({ 2\sinh^2(R) + 1 \over 2 \sinh^2(R) + 3}\biggr)
 \le 4 \alpha^2 x \le {1\over \sinh^2(R) }.$$

Using Proposition (\ref{tuberadbound}) we can, in turn, bound $R$ in terms of $\alpha \l$.
Because of its importance in what follows, we introduce the new variable
$$ \r = h^{-1}(\alpha \l).$$
Note that $\r$ is defined whenever $\alpha \l = h(\r) \leq h_{max}$ and, 
if $R \ge 0.531$ also, it satisfies $0.531 \leq \r \leq R$.  
This allows us to replace $R$ with $\r$ in the estimate (\ref{keybound}):

\begin{proposition}  Whenever $\alpha \l \leq h_{max}$ and 
$R \ge 0.531$ the following inequality holds:
\begin{eqnarray}{-1 \over \sinh^2(\r)}
        \biggl({ 2\sinh^2(\r) + 1 \over 2 \sinh^2(\r) + 3}\biggr)
 \le 4 \alpha^2 x \le {1\over \sinh^2(\r) }.
\label{rhoderl}
\end{eqnarray}
\end{proposition}

\begin{proof}  Proposition \ref{tuberadbound} implies that $\r \leq R$. 
The result follows immediately once it is noted that both
${1\over \sinh^2(r) }$ and ${1\over \sinh^2(r) } 
\biggl({ 2\sinh^2(r) + 1 \over 2 \sinh^2(r) + 3}\biggr)$ are
decreasing in $r$.  That the first is decreasing is obvious; that the
second is decreasing can be seen easily by rewriting it as 
${2 + {1\over \sinh^2(r) }\over 2 \sinh^2(r) + 3}$ so that the numerator
is decreasing and the denominator increasing.
\end{proof}

\medskip

The significance of putting the inequality in the form (\ref{rhoderl}),
as opposed to that of (\ref{keybound}) is that, since
$\r$ is a function of $\alpha \l$, the inequality bounds
the derivative of the core length $\l$ purely in terms of $\alpha$
and $\l$.  Since $\alpha^ 2$ is our parameter, this will allow us
to bound the value $\alpha \l$ by integration, after some algebraic
manipulation and separation of variables.

Now put $$u = \displaystyle{\alpha\over \l}.$$ 
This turns out to be an important and useful function of $\alpha$ and
$\l$.  It approaches a finite, non-zero value as one approaches the 
cusp case, even though $\l$ and $\alpha$ both approach $0$.  
Recall that the meridian length $m$ and annulus height $h$ in the flat
metric on the boundary of a tube of radius $R$ around the singular locus
satisfy $\m =\alpha \sinh R, ~ \h = \l \cosh R$.  Thus, as $R\to \infty$,
the ratio of $\alpha$ to $\l$ approaches that of $\m$ to $\h$.  This
implies that:
\begin{eqnarray}\lim_{R \to \infty} u = \lim_{R \to \infty} {\m \over \h} 
=\lim_{R \to \infty} {\m^2 \over A}= \hat \L^2,
\label{ulim} \end{eqnarray}
where $\hat \L$ is the normalized length of the meridian curve on the
torus boundary of the tube around the cusp.

This provides us with an initial condition for $u$ in terms of the normalized
length of the chosen surgery curve.  To control the value of $\alpha \l$,
it suffices to control the value of $u$. The derivative of $u$ can be
computed by:
$${du \over d\alpha} = {1\over \l} - {\alpha\over \l^2} {d\l\over d\alpha}
={\alpha\over \l^2}\left( {\l\over\alpha} - {d\l\over d\alpha}\right)
= {\alpha\over \l^2} \left(-4\alpha^2 x {\l \over \alpha}\right) 
= -{1\over \l}(4\alpha^2 x),$$
or
$${du \over dt} = {1\over 2 \alpha} {du \over d\alpha}
= -{1\over 2 \alpha \l}(4\alpha^2 x),$$
where $t = \alpha^2$. 

Using (\ref{rhoderl}) and the fact that $h(\r) = \alpha \l$ we obtain
upper and lower bounds on the derivative of $u$ in terms of $\r$.
The expressions for these bounds become simpler if use the variable:
$$ z = \tanh \r.$$
Then, as derived in the proof of Lemma \ref{hinverse}, $h(\r) = 
3.3957 {z (1-z^2) \over 1+z^2}.$ We define the function
\begin{eqnarray}H(z) = {1 \over \alpha \l} = {1 \over h(\r)} = 
{1+z^2 \over {3.3957 z (1-z^2)}}.
\label{Hdef}\end{eqnarray}
Noting that $\sinh^2 (\r) = {z^2 \over 1-z^2}$ we can rewrite the inequality
(\ref{rhoderl}) in terms of $z$:
\begin{eqnarray}-\biggl({(1-z^2)(1+z^2) \over z^2(3-z^2)}\biggr)
\le 4 \alpha^2 x \le {1-z^2 \over z^2}. \label{zderl} \end{eqnarray}
We introduce the functions:
\begin{eqnarray}G(z) = {H(z)\over 2}{1-z^2 \over z^2} = {1+z^2 \over 6.7914 ~ z^3}
\label{Gdef} \end{eqnarray}
and
\begin{eqnarray}\tilde G(z) = {H(z)\over 2}{(1-z^2)(1+z^2) \over z^2(3-z^2)} = 
{(1+z^2)^2 \over 6.7914 ~ z^3\,(3-z^2)}.\label{GGdef}\end{eqnarray} 

Since ${du \over dt} = -{1\over 2 \alpha \l}(4\alpha^2 x) = 
-{H(z) \over 2}(4\alpha^2 x)$, the inequality (\ref{zderl}) provides 
inequalities for ${du \over dt}$ expressed purely in terms of $z$.  Using the 
functions defined above, the inequalities can be written simply as
$$-G(z) \le {du \over dt} \le \tilde G(z).$$

These inequalities hold as long as the tube radius $R\geq 0.531$ 
and $\alpha \l < h_{max}$.  The latter holds, by definition of $h^{-1}$, 
as long as $h^{-1}(\alpha \l) = \r \ge \r_1 = 0.531$, or, 
since $z = \tanh (\r)$ is increasing in $\r$, 
as long as $z \ge \tanh \r_1 \approx 0.4862 = z_1.$
If the initial tube radius is at least $0.531$ then it will remain
so as long as $\r \ge \r_1$.  Thus, in this case, as long as
$z \geq z_1$, the inequalities are valid.  We record this fact as a 
proposition.

\begin{proposition} For any smooth family of hyperbolic cone-manifolds whose initial
tube radius is at least $0.531$, the following differential inequalities
hold as long as $z \geq z_1 = 0.4862$:
\begin{eqnarray} -G(z) \le {du \over dt} \le \tilde G(z),
\label{dudtbound}\end{eqnarray}
where the functions $G(z)$ and $\tilde G(z)$ are
defined by (\ref{Gdef}) and (\ref{GGdef}), respectively.
\end{proposition}

%\bigskip

We will only use the lower bound in this section.  The upper bound 
will be used in the final section.

\medskip

In order to solve this differential inequality, we note that
$u = {\alpha \over \l} = {t \over \alpha \l}$, where $t = \alpha^2$ is
our variational parameter.  By definition, $H(z) = {1\over \alpha \l}$
and this becomes
\begin{eqnarray} u = t H(z)\label{u_z}\end{eqnarray}
From the inequality (\ref{dudtbound}) we obtain
$${d \over dt} ( H(z) t) \ge -G(z)$$
or
\begin{eqnarray}
t {dH \over dz}{dz \over dt} \ge -(H(z) + G(z)).
\label{dzdteq}\end{eqnarray} 

Denoting ${dH \over dz}$ by $H'(z)$ this provides the inequality:
\begin{eqnarray}
{dz \over dt} \ge {-(H(z) + G(z))\over t ~ H'(z)}
\end{eqnarray}

Again, if the initial structure has tube radius at least $0.531$, 
this inequality is valid as long as $z > z_1$.
Observe that $H'(z)$ is {\it positive} since, by Lemma \ref{hinverse},
$h(\r) = {1 \over H(z)}$ is decreasing for these values of $z = \tanh \r$. 

Since this inequality bounds the change in $z$, if we start with a 
complete structure, where $z = 1$, it should provide conditions
under which this inequality will be maintained until $t = (2 \pi)^2$.
In particular, we will have $z\geq z_1$, hence $\alpha \l < h_{max}$,
throughout the deformation, implying, by Theorem \ref{alnondeg}, 
that the smooth structure with cone angle $2 \pi$ can be reached 
without any degeneration. To do this explicitly we will use 
separation of variables.

By algebraic manipulation we obtain 
\begin{eqnarray}
{H'(z) \over H(z) + G(z)} ~ {dz \over dt} \ge -{1 \over t}.
\label{lowerdzbound}\end{eqnarray}
However this separation of
variables is only valid away from the complete structure because both
sides of the new inequality blow up as $t\to 0$ and $z\to 1$.
It cannot be applied directly for initial conditions at the complete
structure.  Some care must be taken to analyze the rate at which the left side
goes to infinity as $t\to 0$.

We compute that $${H'(z) \over H(z) + G(z)} = F(z) + {1\over 1-z}$$
where $$F(z) = -{(1+4z+6z^2+z^4)\over(z+1)(1+z^2)^2},$$
and $F$ is integrable on the interval $0 \le z \le 1$. 
Recall that $z(t)$ is a smooth function of $t$ which approaches
$1$ as $t$ approaches $0$.  For any sufficiently small value of $t>0$,
$z(t)<1$ will be larger than $z_1 = .4862$ and the differential inequality
(\ref{lowerdzbound}) holds.  Choose $0< t_0 < \T$ so that 
$z_1<z(t)<1$ for all $0<t<\T$, and denote $z(t_0)$ by $z_0$.
Integrating the inequality over the interval $0<t<\T$
and changing variable to $w = z(t)$, we obtain
$$\int_{z_0}^{z(\T)} F(w) \,dw ~+~ 
\log(1-z_0)-\log(1-z(\T))\ge \log(t_0) -\log(\T) $$
or 
$$\exp\left(\int_{z_0}^{z(\T)} F(w) \,dw \right) \ge {t_0\over 1-z_0} {1-z(\T)\over \T}.$$

To compute the limit of ${t_0\over 1-z_0}$ as $t_0 \to 0$, multiply
the numerator and denominator by $H(z_0)$. Since $u(t) = H(z(t))t$, this becomes
${u(t_0)\over (1-z_0)H(z_0)}$. 
Now as $t_0\to 0$,  $z_0\to 1$
and from the formula (\ref{Hdef}) it is clear that 
$H(z_0)(1-z_0) \to \displaystyle{1 \over 3.3957}$.  From (\ref{ulim})
we know that $\lim_{R \to \infty} u = \hat \L^2$.  Since $R \to
\infty$ as $t \to 0$, it follows that $\lim_{t \to 0} u(t) = \hat \L^2$. 

We conclude that 
$$\exp\left(\int_1^{z(\T)} F(w) \,dw \right) \ge 3.3957 ~ \hat \L^2 ~{1-z(\T)\over \T}.$$
This inequality holds for any time $\tau$ during a deformation 
through cone-manifolds which begins at a complete structure (where $z(0)=1$),
using a surgery curve of normalized length $\hat \L$, as long as
$z(t)$ is larger than $z_1$ throughout the deformation.  It provides
information about the times $t$ at which various values of $z(t)$
can be attained.
In particular, it implies, for any
$z\geq z_1$, the following inequality for the first time $t$ at which 
$z(t) = z$:
\begin{eqnarray} \label{firsttbound}
t \ge 3.3957 ~ \hat \L^2 ~(1-z)~ \exp\left(-\int_1^z F(w) \,dw \right).
\end{eqnarray}

We conclude that  we can increase the cone angle $\alpha$ from $0$ to $2\pi$, 
maintaining $z=\tanh \r \ge z_1 > \tanh(\r_1)$,
hence keeping the tube radius $R \ge \r \ge \r_1 = 0.531$ and 
$\alpha \l \leq h_{max}$, provided 
$$3.3957 ~ \hat \L^2 ~(1-z_1)~ \exp\left(-\int_1^{z_1} F(w) dw \right)~\ge ~ (2 \pi)^2$$
or
$$ \hat \L^2 ~\ge ~ {(2\pi)^2 \over 3.3957 (1-z_1)} ~ \exp\left(\int_1^{z_1} F(w) dw \right)
\approx 56.4696$$
or 
$$\hat \L ~\ge~ \sqrt{56.4696} \approx 7.5146.$$

Thus, we have shown that as long as the normalized Euclidean geodesic 
length $\hat L$ of the surgery curve satisfies this inequality then
there is a lower bound to the tube radius.  This completes the
proof of Theorem \ref{thm1} which we restate here for convenience.

\begin{theorem}  Let $\M$ be a complete, finite volume
hyperbolic $3$-manifold with one cusp and
let $T$ be a horospherical torus which is embedded as
a cross-section to the cusp. Let
$\gamma$ be a simple closed curve on $T$, $\M(\gamma)$ the Dehn
filling with $\gamma$ as meridian.  Let $\M_\alpha(\gamma)$ be
a cone-manifold structure on $\M(\gamma)$ with cone angle $\alpha$
along the core, $\Sigma$, of the added solid torus, obtained by increasing
the angle from the complete structure.
If the normalized length of $\gamma$ on $T$ is at least $7.515$,
then there is a positive lower bound to the tube radius around $\Sigma$
for all $2\pi \ge \alpha \ge 0$.
\end{theorem}

\begin{remark} \label{geom_of_path} The proof of this theorem shows that $\alpha \l \le
h_{max}$ and $R \ge 0.531$ for $0 \le \alpha \le 2\pi$, where $\l$ and $R$ denote the length
and tube radius of the singular geodesic $\Sigma$.  In particular, the core
geodesic of the non-singular hyperbolic structure on $\M(\gamma)$ has 
length $\l \le {h_{max} \over 2\pi} \approx 0.162$ and tube radius $R \ge 0.531$.
\end{remark}

Theorem \ref{thm1}, together with Theorem \ref{thm2}, implies 
our main result, Theorem \ref{thm3}:
(Note also that Theorem \ref{alnondeg} (which depends on Theorem \ref{thm2}) 
together with (\ref{firsttbound}) immediately implies Theorem \ref{thm3}.)

\begin{theorem}  
Let $\M$ be a complete, orientable hyperbolic $3$-manifold with 
one cusp, and let $T$ be a horospherical torus which is embedded as
a cross-section to the cusp of $\M$. Let
$\gamma$ be a simple closed curve on $T$ whose normalized Euclidean 
geodesic length $\hat \L$ is at least $7.515$. Then the closed manifold 
$\M(\gamma)$ obtained by Dehn filling along $\gamma$ is hyperbolic.
\end{theorem}

%\medskip
This result also gives a universal bound on the
number of non-hyperbolic Dehn fillings on a cusped hyperbolic $3$-manifold.

\begin{corollary} Let $\M$ be a complete, orientable hyperbolic $3$-manifold with 
one cusp. Then at most 60 Dehn fillings on $\M$ yield manifolds which admit no
complete hyperbolic metric.
\end{corollary}

\begin{proof} Let $T$ be the Euclidean torus obtained by taking a horospherical
cusp cross section in $\M$ and rescaling the metric so that $T$ has area $A=1$.
Then (as in Agol \cite{agol}) if $\beta,\gamma$ are slopes with $\M_{\beta}, 
\M{\gamma}$ both non-hyperbolic, the distance $\Delta$ between these slopes satisfies
$\Delta \le (7.515)^2 <57$. Hence $\Delta \le 56$ since $\Delta$ is an integer.
Now application of Lemma 8.2 of \cite{agol} with $p=59$ shows that the number of exceptional
slopes is at most $p+1 = 60$. 
\end{proof}

%\medskip
 
The only way that we have used the hypothesis that there is a single 
cusp is that, when applying the area bound from the packing theorem
(Theorem \ref{areathm}), we used the larger bound (\ref{areabound}). 
This followed from the existence of {\em two} disjoint  embedded 
ellipses in the boundary torus coming from projecting translates
of the tube.  When there are multiple cusps, it can be arranged so 
that there will still be two disjoint embedded ellipses on one 
boundary torus but perhaps only one on the remaining boundary tori. 
(See \cite{BH} for this argument.)  Then Theorem \ref{areathm}
simply gives an area bound (\ref{multipleareabound}) that is half
as large for the remaining boundary tori.

This implies that the function corresponding to $h(r)$ 
(see (\ref{humpeq})) on the remaining boundary tori is half as big.  
It follows that the functions
corresponding to $H(z), G(z), \tilde G(z)$ are twice as big.
Everything else remains the same.  The effect on the 
differential inequalities is that the inequality (\ref{dudtbound}) 
is replaced by one in which $G(z)$ and $\tilde G(z)$ are twice 
as large.  However, the key inequality (\ref{dzdteq}) relating the change
in $z = \tanh(\rho)$ to that of $t = \alpha^2$ remains exactly the
same because it involves the ratio of functions, each of which is
twice as large.  The only change in the analysis arising from that
inequality is that the limit as $z \to 1$ of $H(z) (1-z)$ is
twice as large.  In other words, $H(z)(1-z) \to {2 \over 3.3957}$ 
and the coefficient $3.3957$ in inequality (\ref{firsttbound}) is
replaced by ${3.3957 \over 2}$.  
So, to guarantee that angle $\alpha = 2 \pi$ is reached, we need 
to assume that the normalized Euclidean geodesic lengths of
all of the surgery curves satisfy
$$\hat \L ~\ge~ \sqrt{2 ~(56.4696)} \approx 10.6273.$$
Since the first prime number larger than $2 ~(56.4696) \approx 112.939$ 
is $113$, the bound on the number of exceptional slopes per cusp, for
all but one cusp, becomes $113 + 1 = 114$.  The bound for the other
cusp will still be $60$.

Thus, we have proved that

\begin{theorem}\label{multiplecusp}
Let $\M$ be a complete, finite volume orientable hyperbolic $3$-manifold 
with more than one cusp, and let $T_i$ be a horospherical torus which
is embedded as a cross-section to the $i$th cusp of $\M$. Let
$\{\gamma_i\}$ be simple closed curves on the $T_i$ and suppose that, 
for all $i>1$, the normalized Euclidean geodesic length of
$\gamma_i$ on $T_i$ is at least $10.628$ and for $i=1$ it is
at least $7.515$.  Then the
closed manifold $\M(\gamma)$ obtained by Dehn filling along
$\gamma = \{\gamma_i\}$ is hyperbolic.
In particular, there are at most $60$ choices of $\gamma_1$ on the
first cusp and $114$ choices of $\gamma_i$ on the remaining
cusps so that $\M(\gamma)$ can fail to have a hyperbolic metric.
\end{theorem}

\begin{remark}The results in this section provide initial conditions
which guarantee that any particular collection of cone angles,
all at most $2 \pi$, can be realized by hyperbolic cone-manifold
structures on $\M(\gamma)$.  In particular, they imply the existence
of hyperbolic structures on {\em orbifolds} when the singular locus is
a link.  In this case, the cone angles are all of the form
${{2 \pi} \over n}, n \in \Z$.  The conditions on the the normalized
Euclidean geodesic lengths are replaced by the same condition on
$n$ times the length.  Similarly, the results of Section \ref{compare},
concerning volumes and lengths of the singular locus,
will also apply in this case.
\end{remark}

\section{Geometry comparison}\label{compare}

\subsection{Decreasing the cone angle}\label{decrease}

%\medskip
%
It is natural to ask how general this process of constructing a 
closed hyperbolic manifold is.  Can every closed hyperbolic 
$3$-manifold can be obtained by starting with a non-compact, 
finite volume $3$-manifold with one cusp and increasing the cone angle
from $0$ to $2 \pi$?  Specifically, given a simple closed geodesic $\tau$ in
a closed hyperbolic $3$-manifold $N$, can the cone angle be {\it decreased}
from $2 \pi$ (at the smooth structure) back to angle $0$?
There is no topological obstruction to doing this.
It can be shown (see \cite[Theorem 1.2.1]{kojima}, \cite{agol2})
that $N-\tau$ can be given a complete
finite volume metric with pinched negative curvature so it will be irreducible,
atoroidal, and have infinite fundamental group.  In fact, since it is the
interior of a manifold with non-empty boundary, it is Haken, so Thurston's
geometrization theorem for Haken manifolds implies that it can be given a
hyperbolic structure.  The issue is whether or not the hyperbolic structures
on $N$ and on $N-\tau$ can be connected by a family of hyperbolic cone-manifolds.

In this section we apply our techniques to show that, as long as $\tau$
is sufficiently short, with length less than a universal constant 
independent of $N$, then $N$ can be constructed in this manner.  The
cone angle can be decreased back from $2 \pi$ to $0$.

To see why the condition that $\tau$ be short might arise from the 
techniques of the previous section, note that all of the closed 
hyperbolic manifolds constructed in the proof of Theorem \ref{thm3} 
have a short geodesic, which was the singular set throughout the 
deformation through cone-manifolds.  It is short because the control of
the tube radius using the inverse function $h^{-1}(\alpha \l)$ only held
as long as $\alpha \l \leq h_{max} \approx 1.019675$.  When $\alpha = 2 \pi$ this holds if
$\l \leq 0.162.$

In order to show that it is possible to decrease the cone angle back to
$0$, the main step is again to show that the tube radius is bounded below.
By Theorem \ref{alnondeg}  it suffices to show that the initial tube radius $R$
satisfies $R\geq 0.531$, that the volumes remain bounded, and that, 
if $\alpha \l \leq h_{max}$ at the beginning of the deformation, 
it will remain so throughout.

\begin{lemma} \label{alincrease} $\alpha \l$ is an increasing 
function of $\alpha$ provided the tube radius $R$ satisfies 
$R\ge 0.4407$ and $\alpha \leq 2 \pi$.
\end{lemma}

\begin{proof} Using equation (\ref{dldalpha}) and  estimate 
(\ref{keybound}) (the proof of the estimate using \cite{HK1}
requires that $\alpha \leq 2 \pi$, though this is probably 
unnecessary), we have
$${1\over \alpha} {d(\alpha \l) \over d \alpha}
= {d\l \over d\alpha} + {\l\over \alpha}= {\l \over
\alpha} (2 + 4 \alpha^2 x)
\ge {\l \over \alpha}\left(2 - {1\over \sinh^2 R}\left({2\sinh^2 R + 1 \over
2\sinh^2 R + 3}\right)\right)
\ge 0$$
provided ${1\over s^2}({2s^2 + 1 \over 2s^2 + 3})\le 2$ or
$4s^4+ 4 s^2 -1  \ge 0$ where $s=\sinh R$. This holds provided
$(2s^2+1)^2 \ge 2$, i.e. $s^2 \ge {\sqrt{2}-1\over 2}$ or $R \ge 0.4407$.
\end{proof}

\begin{theorem} \label{connect1} Let $M$ be a closed hyperbolic $3$-manifold 
and $\tau$ a simple closed geodesic in $M$ having length $l(\tau) \le
h_{max}/(2\pi)
\approx 0.1623$ and tube radius $R \ge 0.531$.
Then the hyperbolic structure on $M$ can
be deformed to a complete hyperbolic structure on $M-\tau$ by decreasing 
the cone angle along $\tau$ from $2\pi$ to $0$. 
\end{theorem}

\begin{proof} For $\alpha=2\pi$ we have $\alpha \l \le h_{max}$ and 
$R \ge 0.531> 0.4407$. By Lemma \ref{alincrease}, $\alpha \l \le h_{max}$ 
throughout any deformation decreasing the cone angles.
Since the volume is increasing, it is not immediate that the volumes
are bounded above.  However, it is not difficult to show by general 
arguments, for example by using the Gromov norm, that an upper bound
exists. (See, for example, \cite[Prop. 1.3.2]{kojima} or \cite{agol2}.)
For our specific situation, we can appeal to the next section
in which we give explicit bounds
on the change in volume as the cone angle is changed.  In particular,
this provides an upper bound over the family of cone-manifolds.

Hence there can be no degeneration by Theorem \ref{alnondeg}. 
\end{proof}

The previous theorem requires conditions on both the length of the
geodesic and on its tube radius.  As noted previously, for a general
cone-manifold, it is not true that a sufficiently short 
singular locus provides a lower bound on the tube radius.  However,
for smooth hyperbolic manifolds such a lower bound does exist.  This
follows from the Margulis Lemma or the J{\o}rgensen inequality.

An explicit formula giving a lower bound to the tube radius around
sufficiently short closed geodesics in closed
hyperbolic 3-manifolds was derived by Meyerhoff and Zagier
\cite{meyerhoff} and sharpened by  Cao, Gehring, Martin \cite{CGM}
(see also \cite[Theorem 3.2]{gmm}).
By combining this bound with a tube packing estimate, similar to
the estimate $\alpha \l \ge h(r)$ from the previous section,
Gabai-Milley-Meyerhoff obtain an improved bound
on the tube radius of short geodesics (\cite[Theorem 3.1]{gmm}).
In particular, their
formula implies that if $\tau$ is a closed geodesic in
a smooth hyperbolic $3$-manifold
and if its length satisfies $\l(\tau) \le 0.111$,
then it has tube radius $R \ge 0.982\ge 0.531$. Hence, the previous Theorem
applies. The conclusion is

\begin{corollary} \label{connect2} Let $M$ be a closed hyperbolic
3-manifold and $\tau$ a simple closed geodesic in $M$ having length
$\l(\tau) \le 0.111$.  Then the hyperbolic structure on $M$ can be
deformed to a complete hyperbolic structure on
$M-\tau$ by decreasing the cone angle along
$\tau$ from $2\pi$ to zero.
\end{corollary}

Suppose $\tau$ is actually a {\em shortest} closest geodesic
in a closed hyperbolic 3-manifold $M$. Then $\tau$
is a simple closed curve, and the results of
Gabai-Meyerhoff-Thurston \cite{GMT} show that either
$\tau$ has tube radius $R \ge \log(3)/2 > 0.531$
or $\tau$ has length $> 0.831$.  Thus
if $\tau$ has length $\le 0.162$
then the hypotheses of Theorem \ref{connect1} are again
satisfied. This proves

\begin{corollary} \label{connect3}
Let $M$ be a closed hyperbolic manifold
and let $\tau$ be a shortest closed geodesic in
$M$ having length $\l(\tau) \le 0.162$.
Then the hyperbolic structure on $M$ can
be deformed to a complete hyperbolic structure on $M-\tau$ by decreasing
the cone angle along $\tau$ from $2\pi$ to $0$.
\end{corollary}

\subsection{Volume estimates}\label{volumes}

The rigidity theorem of Mostow and Prasad shows that
geometric invariants of finite volume hyperbolic 3-manifolds are
actually topological invariants. Perhaps the most useful such
invariant is the hyperbolic volume. This volume has proved to be
a good way of distinguishing 3-manifolds,
and is a very good measure of the complexity of a manifold.

Thurston and J{\o}rgensen \cite{thnotes} proved that the
set of volumes of complete, finite volume, orientable, hyperbolic 3-manifolds
is a well-ordered, closed subset of $\R$ of order type $\omega^\omega$,
and that there are finitely many manifolds of any given volume.
Thus the volumes can be arranged:
$$0 <  v_1 < v_2 < \ldots < v_\omega < v_{\omega +1}< \ldots <
v_{2\omega} < \ldots < v_{3\omega} <  \ldots < v_{\omega^2} < \ldots .$$
The smallest volume $v_1$ is the volume of a closed hyperbolic 3-manifold,
and the first limit volume $v_\omega$ represents the volume of the
smallest cusped hyperbolic 3-manifold.

In general, the volume $v$ of each cusped hyperbolic 3-manifold $M$
is a limit point: performing Dehn filling on $M$ produces a collection of
closed hyperbolic manifolds converging geometrically
to the cusped manifold, and their volumes converge to $v$ from below.
Thus the decrease in volume during Dehn filling is an
indication of how close the filled manifold is geometrically
to the cusped manifold.

A few of the lowest volumes are now known. Adams \cite{adams} has 
shown the smallest {\em non-orientable} cusped hyperbolic 3-manifold has 
volume $1.01494 \ldots$; this is the volume of the Gieseking manifold, 
a non-orientable
manifold double covered by the figure eight knot complement.
Recently, Cao-Meyerhoff \cite{cao-meyerhoff} showed that the
{\em orientable} cusped hyperbolic 3-manifolds of smallest volume are the
the figure eight knot complement
and another closely related manifold, with volume $v_\omega = 2.02988 \ldots$.

For closed manifolds, much less is known. The best current
estimate for $v_1$ is that $0.32 < v_1 \le 0.9427 \ldots$,
where right hand side represents the volume of the ``Weeks manifold'' obtained
by $(5,-1),(5,2)$ surgery on the Whitehead link. The left hand side
is an estimate obtained by Agol \cite{agol2}, improving earlier results
of Meyerhoff
\cite{meyerhoff}
\cite{Me}, Gabai-Meyerhoff-Thurston
\cite{GMT}, Gehring-Martin \cite{GM1}, \cite{GM2} and Przeworski \cite{prez}.

Since the smallest cusped manifold volume is known but
the smallest closed manifold volume is not known,
we could try to study volumes of closed hyperbolic 3-manifolds
by regarding them as Dehn fillings on cusped manifolds.
Our work in sections \ref{harmdef} and \ref{tuberad}, gives good control on
the change in length of the core geodesic during Dehn
filling. We now show that this leads to good estimates on
the change in hyperbolic volume during Dehn filling.

%\bigskip

\begin{theorem} \label{volchange} Let $\M$ be a cusped
hyperbolic 3-manifold and $\N$ a closed hyperbolic 3-manifold
which can be joined by a smooth family of 
hyperbolic cone-manifolds
with cone angles $0 \le \alpha \le 2\pi$ along a knot $\Sigma$.
Suppose that $\alpha \l \le h_{max} \approx 1.019675$ holds throughout the
deformation, where $\l$ denotes the length of $\Sigma$.
Then the difference in volume 
$$\Delta V = {\rm Volume}(\M) - {\rm Volume}(\N)$$
satisfies
\begin{eqnarray}\label{voldecrease}
 \int_{\hat z}^{1} {H'(z) dz \over 4 H(z)(H(z) + G(z))} 
\le 
\Delta V
\le \int_{\hat z}^{1} {H'(z) dz \over 4 H(z)(H(z) - \tilde G(z))} 
\end{eqnarray}
where
$\hat z = \tanh(\hat \rho)$, $\hat \rho$
is the unique solution of 
$h(\hat \rho) = 2 \pi \hat\l$ with $\hat \rho \ge 0.531$,
and $\hat \l$ is the length of $\Sigma$ in $\N$ (i.e. when $\alpha=2\pi$).
\end{theorem}

\begin{remark} \label{volgraph}
The graph below shows the upper and lower bounds for $\Delta V$
given by this theorem as a function of the core geodesic 
length $\l$ in $\N$, for  $\l<0.162$. The dotted line shows the asymptotic
formula $\Delta V \sim {\pi \over 2} \l$ as $\l \to 0$
of Neumann-Zagier \cite{neumann-zagier}.

\bigskip
\centerline{\epsfbox{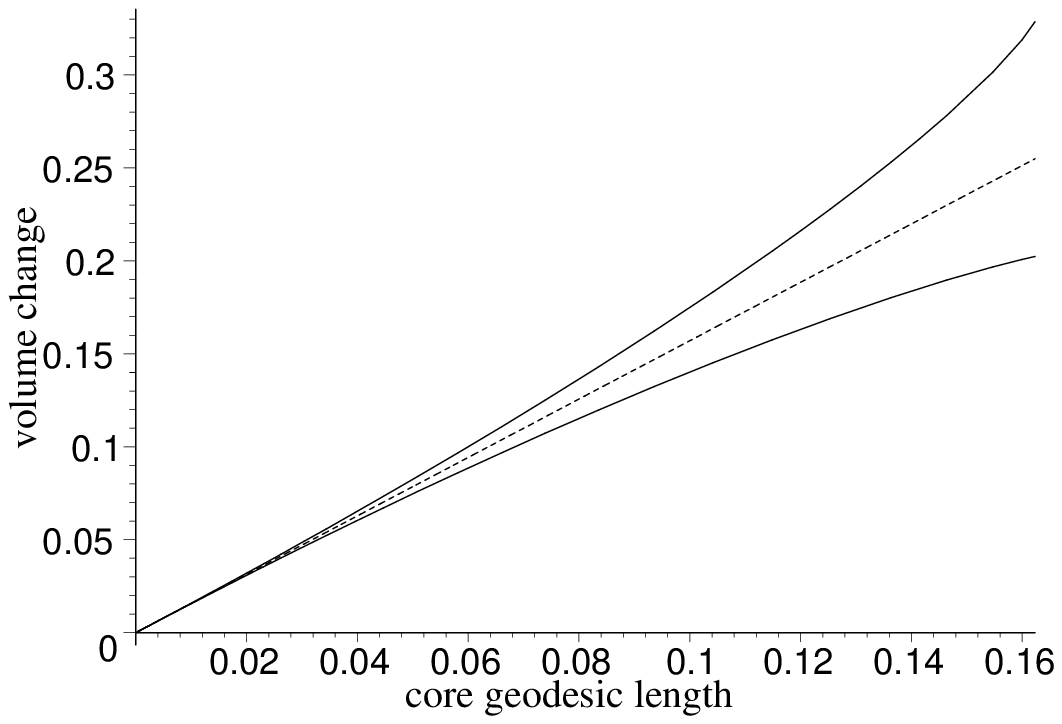}}
%
%\bigskip
\end{remark}

\begin{proof}
We write $t = \alpha^2$ and use the notation from section \ref{tuberad}.
From the Schl\"afli formula (\ref{schlafli}),
the change in volume $V$ of a hyperbolic cone-manifold during a deformation
satisfies
\begin{eqnarray}
dV = - {1\over 2} \l d\alpha = -{\alpha d\alpha \over 2 u} = -{dt \over 4u}
\label{dVol}\end{eqnarray}  
since $\l = {\alpha \over u}$ and $dt = d(\alpha^2) = 2 \alpha d\alpha$. 
Recalling that $u = H(z)t$ we can rewrite this as: 
\begin{eqnarray} {dV \over dt}  = -{1 \over t} {1 \over 4 H(z)}.
\label{HdVol} \end{eqnarray}

Since $\M$ is a cusped manifold, the condition $\alpha\l \le h_{max}$
guarantees that the tube radius satisfies $R \ge \rho_1 = 0.531$ throughout
the deformation (see theorem \ref{alnondeg}).
From equation (\ref{dudtbound}) we have
$$-G(z) \le {du \over dt} \le \tilde G(z),$$
where $G(z)$ and $\tilde G(z)$ are defined as in
(\ref{Gdef}) and (\ref{GGdef}).
Again, since $u = H(z)t$, it follows that  
${du \over dt} = H'(z) t {dz \over dt} + H(z).$
Hence, we obtain:
$$-(G(z)+ H(z)) \le H'(z) t {dz \over dt}  \le \tilde G(z) - H(z).$$

With algebraic manipulation to separate the variables as in 
Section \ref{tuberad} this becomes:
\begin{eqnarray}{H'(z) \over G(z) + H(z)} \,{dz\over dt} \le -{1 \over t} \le
{H'(z) 
\over H(z) - \tilde G(z)} \,{dz\over dt}. \label{dtinequal}\end{eqnarray}
To see that the direction of the inequalities is as claimed,
note as before that $H'(z)$ is positive for all $z > z_1$.  
It is also true, for such values of $z$, that
$H(z) - \tilde G(z)$ is positive.  To see this, recall that
$\tilde G(z) = {H(z)\over 2}~\biggl({(1-z^2)(1+z^2) \over z^2(3-z^2)}\biggr)$.
Hence, it suffices to check that 
${(1-z^2)(1+z^2) \over z^2(3-z^2)}
< 2$.  But ${(1-z^2)(1+z^2) \over z^2(3-z^2)} = 
{1\over \sinh^2 \rho}~
\biggl({2\sinh^2 \rho + 1 \over 2\sinh^2 \rho + 3}\biggr)$ 
and we computed in the proof of 
Lemma \ref{alincrease} that this is less than $2$ as long as
$\rho > 0.4407$.  Since $\rho \ge \rho_1 = 0.531$ the inequality holds.

Putting together equation (\ref{HdVol}) with inequality (\ref{dtinequal}),
we obtain:
\begin{eqnarray}
{H'(z) \over 4H(z)( H(z)+ G(z))}~ {dz\over dt}~ \le ~ {dV\over dt} ~ 
\le ~{H'(z)  \over 4 H(z)(H(z) - \tilde G(z))}~ {dz\over dt}.
\label{volinequal} \end{eqnarray} 

We now integrate over the interval $0 \le t \le \hat t = (2\pi)^2$,
and change variable from $t$ to $z(t)$.  As $t$ increases, the values of 
$z$ {\em decrease} from $z(0)=1$ to 
a value $\hat z = z(\hat t)$, satisfying $\hat z > z_1 = \tanh \rho_1$, so the {\em decrease}
in volume satisfies:
$$
\int_{\hat z}^1 {H'(z) \over 4H(z)( H(z)+ G(z))}~ {dz}~ \le ~ 
\Delta V ~ 
\le \int_{\hat z}^1 {H'(z)  \over 4 H(z)(H(z) - \tilde G(z))}~ {dz}.
$$
(Note that the integrands are positive.)
\end{proof}

The results of section \ref{tuberad} (Remark \ref{geom_of_path}) show that 
Theorem \ref{volchange} applies
when $\N = \M(\gamma)$ is obtained from
a cusped manifold $\M$ by  Dehn filling
along a surgery curve $\gamma$ with  normalized length 
$\hat L \ge 7.515$. This gives

\begin{corollary}\label{volL} Let $\M$ be complete, finite volume 
hyperbolic manifold with one cusp and $\gamma$ a surgery curve
with normalized length $\hat L \ge 7.515$.  Then $\M(\gamma)$
is hyperbolic and its volume satisfies: 
$${\rm Volume}(\M(\gamma)) \ge {\rm Volume}(\M) - 0.329.$$
In particular,
$${\rm Volume}(\M(\gamma)) \ge 1.701.$$
\end{corollary}

\begin{proof}
For surgery curves $\gamma$ with 
normalized length $\hat \L \ge7.515$, we have
$\alpha \l \le h_{max}$ and tube radius $R\ge \r_1 = 0.531$ 
as the cone angle
is increased from $0$ to $2 \pi$.  
The values of $z$ {\em decrease} from $1$ to 
a value $\hat z$, satisfying $\hat z > z_1 = \tanh \r_1$, so the {\em decrease}
in volume is at most:
$$\int_{\hat z}^1 {H'(z) dz \over 4 H(z)(H(z) - \tilde G(z))} 
\leq \int_{z_1}^1 {H'(z) dz \over 4 H(z)(H(z) - \tilde G(z))} < 0.3287,$$
i.e. $${\rm Volume}(\M(\gamma)) \ge {\rm Volume}(\M) - 0.3287.$$

The results of Cao-Meyerhoff \cite{cao-meyerhoff}
show that the figure eight knot complement
and its sister are the cusped orientable hyperbolic $3$-manifolds of minimal volume
$\approx 2.02988$. So we conclude that any surgery on a cusped hyperbolic 
$3$-manifold along a surgery curve with $\hat L \ge 7.515$ gives a hyperbolic manifold with
volume at least
$2.0298 - 0.3287 = 1.701$.
\end{proof}

The results of section \ref{decrease} show that we can also
apply Theorem \ref{volchange} when $M$ is a closed
hyperbolic 3-manifold and $X=M-\tau$ is obtained
by removing a sufficiently short simple closed geodesic $\tau$.
For example, the proof of Corollary \ref{connect3}
shows that the hypotheses of Theorem \ref{volchange}
 are satisfied for any {\em shortest} closed
geodesic $\tau$ of length at most  $0.162$. By the same argument
as above, we then obtain the
following estimate on volumes of closed 3-manifolds
containing short geodesics.

\begin{corollary} \label{volshort}
Let $M$ be a closed hyperbolic $3$-manifold
and let $\tau$ be a shortest closed geodesic in
$M$ having length $\l(\tau) \le 0.162$.
Then $M - \tau$ has a finite volume hyperbolic structure
and
$${\rm Volume}(M) \ge {\rm Volume}(M - \tau) - 0.329.$$
In particular,
$${\rm Volume}(M) \ge 1.701.$$
\end{corollary}

\begin{remark}
It is interesting to compare this with the result of
Agol in \cite{agol2}, which shows that any closed hyperbolic
$3$-manifold with shortest geodesic length less than $0.244$
has volume greater than the volume of the Weeks manifold ($0.9427\ldots$).
That there is such a large gap between the volume estimate above
($1.701$) and the volume of the Weeks manifold suggests that one
ought to be able to improve significantly our bound on the
geodesic length.  Unfortunately, the fact that all our
arguments using the function $h(r)$ break down once the
geodesic length gets much bigger means that such an improvement
would require further methods.
\end{remark}

\begin{remark} We saw in Remark \ref{lincreases} that for a
hyperbolic cone-manifold, if the tube radius $R$ satisfies
$R \geq {\rm arcsinh}({1 \over \sqrt{2}}) \approx 0.6584$,
then the core length $\l$ decreases as $\alpha$ decreases.
If the length of a closed geodesic in a closed
hyperbolic $3$-manifold satisfies
$\l(\tau) \le 0.111$, then $\tau$ has tube radius
$R \ge \rho(\tau)> 0.98> 0.6584$.  Furthermore, we can
decrease the cone angle $\alpha$ along $\tau$ from $2 \pi$ to zero,
keeping the tube radius larger than this value throughout the
deformation.  As discussed in Remark \ref{lincreases}, this implies,
by Schl\"afli's formula, that
$${\rm Volume}(M - \tau) \leq {\rm Volume}(M) + \pi \l(\tau).$$
However, it is not hard to see that the estimate 
given in Theorem \ref{volchange} is considerably stronger.
(Compare the figure in Remark \ref{volgraph}.)

Bridgeman (\cite{bridgeman}) showed that such an estimate
holds for a certain nice class of geodesics in hyperbolic 3-manifolds.
However, the estimate does not hold in general. Using Oliver Goodman's
``Tube" program \cite{tube}, Ian Agol has observed that
that this estimate is violated for several closed geodesics $\tau$ in
the Weeks manifold (see \cite{agol2}).
\end{remark}

\end{document}